\documentclass[12pt]{article}%
\usepackage{amsfonts}
\usepackage{amsmath}
\usepackage{amssymb}
\usepackage{graphicx}%
\newtheorem{theorem}{Theorem}
\newtheorem{acknowledgement}[theorem]{Acknowledgement}

\newtheorem{corollary}[theorem]{Corollary}

\newtheorem{definition}[theorem]{Definition}

\newtheorem{lemma}[theorem]{Lemma}

\newtheorem{proposition}[theorem]{Proposition}
\newtheorem{remark}[theorem]{Remark}

\newenvironment{proof}[1][Proof]{\noindent\textbf{#1.} }{\ \rule{0.5em}{0.5em}}
\begin{document}

\title{{\huge Viability for differential equations driven by fractional Brownian
motion}}
\author{{\large Ioana Ciotir}$^{\text{ }a}$$\quad\quad${\large \medskip\medskip
\medskip Aurel R\u{a}\c{s}canu}$^{\text{ }b}$\\$^{a\text{ }}${\small Department of Economics, "Al. I. Cuza" University, Bd.
Carol no. 9-11, Ia\c{s}i, }\\{\small Rom\^{a}nia, e-mail: ioana.ciotir@feaa.uaic.ro\medskip\ }\\$^{b\text{ }}${\small Department of Mathematics, "Al. I. Cuza" University, Bd.
Carol no. 9-11, Ia\c{s}i, \&}\\{\small "Octav Mayer" Mathematics Institute of the Romanian Academy, Bd. Carol
I, no.8, }\\{\small Romania, e-mail: aurel.rascanu@uaic.ro\medskip} }
\maketitle

\begin{abstract}
In this paper we prove a viability result for multidimensional, time
dependent, stochastic differential equations driven by fractional Brownian
motion with Hurst parameter$\frac{1}{2}<$ $H<1,$ using pathwise approach. The
sufficient condition is also an alternative global existence result for the
fractional differential equations with restrictions on the state.

\end{abstract}

\begin{acknowledgement}
{\small The work for this paper was supported by the founds from the Grant
ID\_395/2007 and Grant CEEX, contract CERES-2-Cex 06-11-56/2006.}
\end{acknowledgement}

\medskip


\section{Introduction}

Let $B=\{B_{t},t\geq0\}$ be a fractional Brownian motion (fBm) of Hurst
parameter $H\in(0,1)$. That is, $B$ is a centered Gaussian process with the
covariance function (see \cite{MaVafr68})
\begin{equation}
R_{H}(s,t)=\mathbb{E}\left(  B_{s}B_{t}\right)  =\frac{1}{2}\left(
t^{2H}+s^{2H}-\left\vert t-s\right\vert ^{2H}\right)  . \label{covari}%
\end{equation}
Notice that if $H=\frac{1}{2}$, the process $B$ is a standard Brownian motion,
but if $H\neq\frac{1}{2}$, it\ does not have independent increments. It
clearly follows that $\mathbb{E}|B_{t}-B_{s}|^{2}=|t-s|^{2H}$. As a
consequence, the process $B$ has $\alpha-$Holder continuous paths for all
$\alpha\in(0,H).$

The definition of stochastic integrals with respect to the fractional Brownian
motion has been investigated by several authors\bigskip.

Essentially two different types of integrals can be defined:

\begin{itemize}
\item The divergence integral (or Skorohod integral) with respect to fBm is
defined as the adjoint of the derivative operator in the framework of the
Malliavin calculus. This approach was introduced by Decreusefond and
\"{U}st\"{u}nel \cite{du} and developed by Carmona and Coutin \cite{cc},
Duncan, Hu and Pasik-Duncan \cite{dhp}, Alos, Mazet and Nualart \cite{amn}, Hu
and \O ksendal \cite{HO}, among others. This stochastic integral can be
expressed as the limit of Riemann sums defined using Wick products and
satisfies the zero mean property

\item The pathwise Riemann-Stieltjes integral $%
%TCIMACRO{\tint _{0}^{T}}%
%BeginExpansion
{\textstyle\int_{0}^{T}}
%EndExpansion
u_{s}dB_{s}^{H}$ which exists if the stochastic process $\left(  u_{t}\right)
_{t\in\left[  0,T\right]  }$ has continuous paths of order $\alpha>1-H$, is a
consequence of the result of Young \cite{y}. Z\"{a}hle has defined in
\cite{Za98} and \cite{Za-II} this integral for processes with paths in a
fractional Sobolev type space. In this paper we will follow this last approach.
\end{itemize}

The aim of this paper is to state necessary and sufficient conditions that
guarantee that the solution of a given (forward) stochastic differential
equation \ driven by the fractional Brownian motion $B$ with Hurst parameter
$1/2<H<1$ (in short: f-SDE), $\mathbb{P}-a.s.\;\omega\in\Omega$%
\begin{equation}
X_{s}^{t,x}=x+\int_{t}^{s\vee t}b(r,X_{r}^{t,x})dr+\int_{t}^{s\vee t}%
\sigma(r,X_{r}^{t,x})dB_{r}^{H},\quad s\in\lbrack t,T], \label{fBm}%
\end{equation}
$\left(  t,x\right)  \in\lbrack0,T]\times\mathbb{R}^{k}$ evolves in a
prescribed set $K$ i.e., under which it holds that for all $t\in\lbrack0,T]$
and for all $x\in K:$
\[
X_{s}^{t,x}\in K\quad a.s.\;\omega\in\Omega,\;\forall s\in\lbrack t,T].
\]
where

\begin{itemize}
\item $B=\left(  B^{i}\right)  _{k\times1}$ , $B^{i}$, $i=\overline{1,k},$ are
independent fractional Brownian motions with Hurst parameter $H,$ $\frac{1}%
{2}<H<1.$, and the integral with respect to $B$ is a pathwise
Riemann-Stieltjes integral;

\item $X_{0}=\left(  X_{0}^{i}\right)  _{d\times1}$ is a $d$ - dimensional
random variable defined in a complete probability space $(\Omega
,\mathcal{F},\mathbb{P)};$

\item $b\left(  \omega,\cdot,\cdot\right)  :\left[  0,T\right]  \times
\mathbb{R}^{d}\mathbb{\rightarrow}\mathbb{R}^{d}$ and $\sigma\left(
\omega,\cdot,\cdot\right)  :\left[  0,T\right]  \times\mathbb{R}%
^{d}\mathbb{\rightarrow R}^{d\times k}$ are continuous functions.
\end{itemize}

A general result on the existence and uniqueness of the solution for
multidimensional, time dependent, stochastic differential equations driven by
a fractional Brownian motion with Hurst parameter $H>1/2$ has been given by
Nualart, R\u{a}\c{s}canu in \cite{nr} using a techniques of the classical
fractional calculus.

The viability - property, has been extensively studied for deterministic
differential equations and inclusions,starting with Nagumo's pioneer work in
1943 ( see \cite{av} for references). To our knowledge the first work that
gives a characterization of the viability property in a stochastic framework
was written by Aubin and Da Prato \cite{ADP1} in 1990 (see also \cite{ADP2},
\cite{GT}). The key point of their work consists in defining a suitable
\emph{Bouligand's stochastic tangent cone} which generalizes the cone used
\ in the study of the viability property for deterministic systems. (for other
points of view see also \cite{Mi}, \cite{Mich} and \cite{Mot}).

Another approach has been developed by Buckdahn, Quincampoix, Rainer,
R\u{a}\c{s}canu in \cite{bqrr} which present a unified approach for the study
of the viability property of SDE, BSDE and PDE by relating the distance to the
constraint to some suitable PDE. More precisely, in the case of an SDE
(respectively BSDE) they show that the viability property is equivalent to the
fact that the square of the distance function is a viscosity super solution
(respectively subsolution) of the PDE.

In this paper we will prove a Nagumo type Theorem on viability proprieties of
close bounded subsets with respect to a stochastic differential equation
driven by fractional Brownian motion, following an approach inspired by the
work of Nualart and\ R\u{a}\c{s}canu \cite{nr}.

The organization of the paper is as follows. In Section 2 we recall some
classical definitions and consider the assumptions on the coefficients
supposed to hold. In Section 3 we state our main result. Section 4 contains
the deterministic more general result and Section 5 contains the proof of the
result in the stochastic case.

\section{Preliminaries}

\subsection{Generalized Stieltjes integral}

Let $d,k\in\mathbb{N}^{\ast}.$Given a matrix $A=\left(  a^{i,j}\right)
_{d\times k}$ and a vector $y=\left(  y^{i}\right)  _{d\times1}$ we denote
$\left\vert A\right\vert ^{2}=\sum_{i,j}\left\vert a^{i,j}\right\vert ^{2}$
and $\left\vert y\right\vert ^{2}=\sum_{i}\left\vert y^{i}\right\vert ^{2}.$

Let $t\in\left[  0,T\right]  $ be fixed. Denote by $W^{\alpha,\infty
}(t,T;\mathbb{R}^{d}),$ $0<\alpha<1,$ the space of continuous functions
$f:[t,T]\rightarrow\mathbb{R}^{d}$ such that%
\[
\left\Vert f\right\Vert _{\alpha,\infty;\left[  t,T\right]  }:=\sup
_{s\in\lbrack t,T]}\left(  |f(s)|+\int_{t}^{s}\dfrac{\left\vert f\left(
s\right)  -f\left(  r\right)  \right\vert }{\left(  s-r\right)  ^{\alpha+1}%
}dr\right)  <\infty.
\]
A equivalent norm can be defined by%
\[
\left\vert \!\left\vert \!\left\vert f\right\vert \!\right\vert \!\right\vert
_{\alpha,\lambda;\left[  t,T\right]  }\overset{def}{=}\sup_{s\in\lbrack
t,T]}e^{-\lambda s}\left(  |f(s)|+%
%TCIMACRO{\dint _{t}^{s}}%
%BeginExpansion
{\displaystyle\int_{t}^{s}}
%EndExpansion
\dfrac{\left\vert f\left(  s\right)  -f\left(  r\right)  \right\vert }{\left(
s-r\right)  ^{\alpha+1}}dr\right)
\]
for any $\lambda\geq0.$\newline For any $0<\mu\leq1$, denote by $C^{\mu
}\left(  \left[  t,T\right]  ;\mathbb{R}^{d}\right)  $ the space of $\mu$-
Holder continuous functions $f:[t,T]\rightarrow\mathbb{R}^{d}$, equipped with
the norm%
\[
\left\Vert f\right\Vert _{\mu;\left[  t,T\right]  }\overset{def}{=}\left\Vert
f\right\Vert _{\infty;\left[  t,T\right]  }+\sup_{t\leq s<r\leq T}%
\dfrac{\left\vert f\left(  s\right)  -f\left(  r\right)  \right\vert }{\left(
s-r\right)  ^{\mu}}<\infty,
\]
where $\left\Vert f\right\Vert _{\infty;\left[  t,T\right]  }:=\sup
\limits_{s\in\lbrack t,T]}|f(s)|$. We have, for all $0<\varepsilon<\alpha$%
\[
C^{\alpha+\varepsilon}(\left[  t,T\right]  ;\mathbb{R}^{d})\subset
W^{\alpha,\infty}(t,T;\mathbb{R}^{d})\subset C^{\alpha-\varepsilon}(\left[
t,T\right]  ;\mathbb{R}^{d}).
\]
with continuous embeddings.

Fix a parameter $0<\alpha<\frac{1}{2}$. Denote by $\tilde{W}^{1-\alpha,\infty
}(t,T;\mathbb{R}^{d})$ the space of continuous functions $g:[t,T]\rightarrow
\mathbb{R}^{k}$ such that%
\[
\left\Vert g\right\Vert _{\tilde{W}^{1-\alpha,\infty}(t,T;\mathbb{R}^{d}%
)}:=\left\vert g\left(  t\right)  \right\vert +\sup_{t<r<s<T}\left(
\frac{|g(s)-g(r)|}{(s-r)^{1-\alpha}}+\int_{r}^{s}\frac{|g(y)-g(r)|}%
{(y-r)^{2-\alpha}}dy\right)  <\infty.
\]
Clearly,%
\[
C^{1-\alpha+\varepsilon}(\left[  t,T\right]  ;\mathbb{R}^{d})\subset\tilde
{W}^{1-\alpha,\infty}(t,T;\mathbb{R}^{d})\subset C^{1-\alpha}(\left[
t,T\right]  ;\mathbb{R}^{d}),\quad\forall\varepsilon>0.
\]
Denoting%
\[
\Lambda_{\alpha}(g;\left[  t,T\right]  )\overset{def}{=}\frac{1}%
{\Gamma(1-\alpha)}\sup_{t<r<s<T}\left\vert \left(  D_{s-}^{1-\alpha}%
g_{s-}\right)  (r)\right\vert ,
\]
where%
\[%
\begin{array}
[c]{l}%
\Gamma(\alpha)=\int_{0}^{\infty}p^{\alpha-1}e^{-p}dp\text{ is the Euler
function and}\medskip\\
\left(  D_{s-}^{1-\alpha}g_{s-}\right)  \left(  r\right)  =\dfrac
{e^{i\pi\left(  1-\alpha\right)  }}{\Gamma\left(  \alpha\right)  }\left(
\dfrac{g\left(  r\right)  -g\left(  s\right)  }{\left(  s-r\right)
^{1-\alpha}}+\left(  1-\alpha\right)
%TCIMACRO{\dint _{r}^{s}}%
%BeginExpansion
{\displaystyle\int_{r}^{s}}
%EndExpansion
\dfrac{g\left(  r\right)  -g\left(  y\right)  }{\left(  y-r\right)
^{2-\alpha}}dy\right)  1_{\left(  t,s\right)  }(r)
\end{array}
\]
we have
\[
\Lambda_{\alpha}(g;\left[  t,T\right]  )\leq\frac{1}{\Gamma(1-\alpha
)\Gamma(\alpha)}\left\Vert g\right\Vert _{\tilde{W}^{1-\alpha,\infty
}(t,T;\mathbb{R}^{d})}.
\]
Note that%
\[
\Lambda_{\alpha}(g;\left[  t,T\right]  )\leq\Lambda_{\alpha}(g;\left[
0,T\right]  )\left(  \overset{def}{=}\Lambda_{\alpha}(g)\right)  .
\]
Let $W^{\alpha,1}(t,T;\mathbb{R}^{d})$ the space of measurable functions $f$
on $[t,T]$ such that
\[
\left\Vert f\right\Vert _{\alpha,1;\left[  t,T\right]  }\overset{def}{=}%
\int_{t}^{T}\left[  \frac{|f(s)|}{\left(  s-t\right)  ^{\alpha}}+\int_{t}%
^{s}\frac{|f(s)-f(y)|}{(s-y)^{\alpha+1}}dy\right]  ds<\infty.
\]
Clearly $W^{\alpha,\infty}(t,T;\mathbb{R}^{d})\subset W^{\alpha,1}%
(t,T;\mathbb{R}^{d})$ and $\left\Vert f\right\Vert _{\alpha,1;\left[
t,T\right]  }\leq\left(  T+\frac{T^{1-\alpha}}{1-\alpha}\right)  \left\Vert
f\right\Vert _{\alpha,\infty;\left[  t,T\right]  }~.$ Denote%
\begin{equation}
\left(  D_{t+}^{\alpha}f\right)  \left(  r\right)  =\dfrac{1}{\Gamma\left(
1-\alpha\right)  }\left(  \dfrac{f\left(  r\right)  }{\left(  r-t\right)
^{\alpha}}+\alpha\int_{t}^{r}\dfrac{f\left(  r\right)  -f\left(  y\right)
}{\left(  r-y\right)  ^{\alpha+1}}dy\right)  1_{\left(  t,T\right)  }(r).
\label{z9}%
\end{equation}

\begin{definition}
\label{def-int}Let $0<\alpha<\frac{1}{2}.$ If $f\in W^{\alpha,1}%
(t,T;\mathbb{R}^{dxk})$ and $g\in\tilde{W}^{1-\alpha,\infty}(t,T;\mathbb{R}%
^{k})$, then defining
\begin{equation}
\int_{t}^{s}f\left(  r\right)  dg\left(  r\right)  \overset{def}{=}\left(
-1\right)  ^{\alpha}\int_{t}^{s}\left(  D_{t+}^{\alpha}f\right)  (r)\left(
D_{s-}^{1-\alpha}g_{s-}\right)  \left(  r\right)  dr. \label{integrala}%
\end{equation}
the integral $%
%TCIMACRO{\dint _{t}^{s}}%
%BeginExpansion
{\displaystyle\int_{t}^{s}}
%EndExpansion
fdg$ exists for all $s\in\lbrack t,T]$ and%
\begin{equation}%
\begin{array}
[c]{lll}%
\left\vert
%TCIMACRO{\dint _{t}^{T}}%
%BeginExpansion
{\displaystyle\int_{t}^{T}}
%EndExpansion
f\left(  r\right)  dg\left(  r\right)  \right\vert  & \leq & \sup
\limits_{t\leq r<s\leq T}\left\vert \left(  D_{s-}^{1-\alpha}g_{s-}\right)
(r)\right\vert
%TCIMACRO{\dint _{t}^{T}}%
%BeginExpansion
{\displaystyle\int_{t}^{T}}
%EndExpansion
\left\vert \left(  D_{t+}^{\alpha}f\right)  (s)\right\vert \ ds\medskip\\
& \leq & \Lambda_{\alpha}(g;\left[  t,T\right]  )\left\Vert f\right\Vert
_{\alpha,1;\left[  t,T\right]  }~.
\end{array}
\label{est-int}%
\end{equation}

\end{definition}

\bigskip

We give two continuity properties of this integral used in the paper. Let
$t\leq s<\tau\leq T.$ We denote$\medskip$%
\[
G_{s,\tau}\left(  f\right)  =\int_{s}^{\tau}f\left(  r\right)  dg\left(
r\right)  =G_{t,\tau}\left(  f\right)  -G_{t,s}\left(  f\right)
\]

\begin{proposition}
\label{p-est}Let $0<\alpha<\dfrac{1}{2}$.\newline$\left(  a\right)  \quad
G_{t,T}:W^{\alpha,\infty}(t,T;\mathbb{R}^{d})\rightarrow C^{1-\alpha}(\left[
t,T\right]  ;\mathbb{R}^{d})$ is a linear continuous map and%
\begin{equation}
\left\Vert G_{t,\cdot}\left(  f\right)  \right\Vert _{1-\alpha;\left[
t,T\right]  }\leq A_{\alpha,T}^{\left(  1\right)  }~\Lambda_{\alpha}(g;\left[
t,T\right]  )\left\Vert f\right\Vert _{\alpha,\infty;\left[  t,T\right]  }
\label{est10}%
\end{equation}
where $A_{\alpha,T}^{\left(  1\right)  }$ is a positive constant depending
only on $\alpha$ and $T$; $A_{\alpha,T}^{\left(  1\right)  }\leq4+3T$%
.\newline$\left(  b\right)  \quad G_{t,T}:W^{\alpha,\infty}(t,T;\mathbb{R}%
^{d})\rightarrow W^{\alpha,\infty}(t,T;\mathbb{R}^{d})$ is a linear continuous
map and for all $\lambda>1$%
\begin{equation}
\left\vert \!\left\vert \!\left\vert G_{t,\cdot}\left(  f\right)  \right\vert
\!\right\vert \!\right\vert _{\alpha,\lambda;\left[  t,T\right]  }\leq
\dfrac{\Lambda_{\alpha}(g;\left[  t,T\right]  )}{\lambda^{1-2\alpha}}%
A_{\alpha,T}^{\left(  2\right)  }~\left\Vert f\right\Vert _{\alpha
,\lambda;\left[  t,T\right]  }, \label{est11}%
\end{equation}
where $A_{\alpha,T}^{\left(  2\right)  }$ is a positive constant depending
only on $\alpha$ and $T$; $A_{\alpha,T}^{\left(  2\right)  }=\dfrac
{4}{1-2\alpha}\left(  \dfrac{2}{\alpha}+T^{\alpha}\right)  .$
\end{proposition}

\begin{proof}
$\left(  a\right)  \quad$Let $t\leq\tau<s\leq T.$ From the definition of the
integral we have%
\begin{equation}
\left\vert G_{\tau,s}\left(  f\right)  \right\vert \leq\Lambda_{\alpha
}(g;\left[  t,T\right]  )%
%TCIMACRO{\dint _{\tau}^{s}}%
%BeginExpansion
{\displaystyle\int_{\tau}^{s}}
%EndExpansion
\left(  \dfrac{|f(\theta)|}{(\theta-\tau)^{\alpha}}+\alpha%
%TCIMACRO{\dint _{\tau}^{\theta}}%
%BeginExpansion
{\displaystyle\int_{\tau}^{\theta}}
%EndExpansion
\dfrac{|f(\theta)-f(u)|}{(\theta-u)^{\alpha+1}}du\right)  d\theta\label{est1}%
\end{equation}
and therefore%
\begin{equation}
\left\vert G_{\tau,s}\left(  f\right)  \right\vert \leq\Lambda_{\alpha
}(g;\left[  t,T\right]  )\left(  2+T^{\alpha}\right)  \left(  s-\tau\right)
^{1-\alpha}\left\Vert f\right\Vert _{\alpha,\infty;\left[  t,T\right]  }~.
\label{est2}%
\end{equation}
Hence the inequality (\ref{est10}) follows with $A_{\alpha,T}^{\left(
1\right)  }=\frac{T^{1-\alpha}}{1-\alpha}+T+2+T^{\alpha}\leq4+3T.$

$\left(  b\right)  $ $\quad$We have\medskip\newline$%
%TCIMACRO{\dint _{t}^{s}}%
%BeginExpansion
{\displaystyle\int_{t}^{s}}
%EndExpansion
\dfrac{\left\vert G_{t,s}\left(  f\right)  -G_{t,r}\left(  f\right)
\right\vert }{\left(  s-r\right)  ^{1+\alpha}}dr$\medskip\newline$%
\begin{array}
[c]{c}%
\quad\quad
\end{array}
\leq\Lambda_{\alpha}(g;\left[  t,T\right]  )%
%TCIMACRO{\dint _{t}^{s}}%
%BeginExpansion
{\displaystyle\int_{t}^{s}}
%EndExpansion
\left(  s-r\right)  ^{-\alpha-1}\times\left(
%TCIMACRO{\dint _{r}^{s}}%
%BeginExpansion
{\displaystyle\int_{r}^{s}}
%EndExpansion
\dfrac{|f(\theta)|}{(\theta-r)^{\alpha}}d\theta+\alpha%
%TCIMACRO{\dint _{r}^{\theta}}%
%BeginExpansion
{\displaystyle\int_{r}^{\theta}}
%EndExpansion
\dfrac{|f(\theta)-f(u)|}{(\theta-u)^{\alpha+1}}dud\theta\right)  dr$%
\medskip\newline$%
\begin{array}
[c]{c}%
\quad\quad
\end{array}
\leq\Lambda_{\alpha}(g;\left[  t,T\right]  )%
%TCIMACRO{\dint _{t}^{s}}%
%BeginExpansion
{\displaystyle\int_{t}^{s}}
%EndExpansion
|f(\theta)|%
%TCIMACRO{\dint _{t}^{\theta}}%
%BeginExpansion
{\displaystyle\int_{t}^{\theta}}
%EndExpansion
\left(  s-r\right)  ^{-\alpha-1}(\theta-r)^{-\alpha}drd\theta+$\medskip
\newline$%
\begin{array}
[c]{c}%
\quad\quad
\end{array}
+\Lambda_{\alpha}(g)\alpha%
%TCIMACRO{\dint _{t}^{s}}%
%BeginExpansion
{\displaystyle\int_{t}^{s}}
%EndExpansion%
%TCIMACRO{\dint _{t}^{\theta}}%
%BeginExpansion
{\displaystyle\int_{t}^{\theta}}
%EndExpansion
\dfrac{|f(\theta)-f(u)|}{(\theta-u)^{\alpha+1}}\left(
%TCIMACRO{\dint _{t}^{u}}%
%BeginExpansion
{\displaystyle\int_{t}^{u}}
%EndExpansion
\left(  s-r\right)  ^{-\alpha-1}dr\right)  dud\theta.$\medskip\newline Since
\begin{align*}
&  \int_{t}^{\theta}\left(  s-r\right)  ^{-\alpha-1}(\theta-r)^{-\alpha}dr\\
&  =\left(  s-\theta\right)  ^{-2\alpha}\int_{0}^{\frac{\theta-t}{s-\theta}%
}\left(  1+u\right)  ^{-\alpha-1}u^{-\alpha}du\\
&  \leq\left(  s-\theta\right)  ^{-2\alpha}\left(  \int_{0}^{1}\left(
1+u\right)  ^{-\alpha-1}u^{-\alpha}du+\int_{1}^{\infty}\left(  1+u\right)
^{-\alpha-1}u^{-\alpha}du\right)  \\
&  \leq\left(  \frac{1}{1-\alpha}+\frac{1}{\alpha}\right)  \left(
s-\theta\right)  ^{-2\alpha}.
\end{align*}
and for $t<r<u<\theta<s$%
\[
\int_{t}^{u}\left(  s-r\right)  ^{-\alpha-1}dr\leq\frac{1}{\alpha}\left(
s-u\right)  ^{-\alpha}\leq\frac{T^{\alpha}}{\alpha}\left(  s-t\right)
^{-2\alpha},
\]
then it follows\medskip\newline$\left\vert G_{t,s}\left(  f\right)
\right\vert +%
%TCIMACRO{\dint _{t}^{s}}%
%BeginExpansion
{\displaystyle\int_{t}^{s}}
%EndExpansion
\dfrac{\left\vert G_{t,s}\left(  f\right)  -G_{t,r}\left(  f\right)
\right\vert }{\left(  s-r\right)  ^{1+\alpha}}dr\medskip\newline%
\begin{array}
[c]{c}%
\quad\quad
\end{array}
\leq\Lambda_{\alpha}\left(  g;\left[  t,T\right]  \right)  \left(  \dfrac
{1}{\left(  1-\alpha\right)  \alpha}+T^{\alpha}\right)  $\medskip\newline$%
\begin{array}
[c]{c}%
\quad\quad\quad
\end{array}
\times%
%TCIMACRO{\dint _{t}^{s}}%
%BeginExpansion
{\displaystyle\int_{t}^{s}}
%EndExpansion
\left(  \left(  s-\theta\right)  ^{-2\alpha}+\left(  \theta-t\right)
^{-\alpha}\right)  \left(  \left\vert f(\theta)\right\vert +%
%TCIMACRO{\dint _{t}^{\theta}}%
%BeginExpansion
{\displaystyle\int_{t}^{\theta}}
%EndExpansion
\dfrac{|f(\theta)-f(u)|}{(\theta-u)^{\alpha+1}}du\right)  d\theta.$%
\medskip\newline The inequality (\ref{est11}) clearly follows since$\medskip
\newline%
%TCIMACRO{\dint _{t}^{s}}%
%BeginExpansion
{\displaystyle\int_{t}^{s}}
%EndExpansion
e^{-\lambda\left(  s-\theta\right)  }\left[  \left(  s-\theta\right)
^{-2\alpha}+\left(  \theta-t\right)  ^{-\alpha}\right]  d\theta\medskip
\newline%
\begin{array}
[c]{c}%
\quad\quad\quad\quad
\end{array}
=\dfrac{1}{\lambda^{1-2\alpha}}%
%TCIMACRO{\dint _{0}^{\lambda\left(  s-t\right)  }}%
%BeginExpansion
{\displaystyle\int_{0}^{\lambda\left(  s-t\right)  }}
%EndExpansion
e^{-u}u^{-2\alpha}du+\dfrac{1}{\lambda^{1-\alpha}}e^{-\lambda\left(
s-t\right)  }%
%TCIMACRO{\dint _{0}^{\lambda\left(  s-t\right)  }}%
%BeginExpansion
{\displaystyle\int_{0}^{\lambda\left(  s-t\right)  }}
%EndExpansion
e^{u}u^{-\alpha}du\medskip\newline%
\begin{array}
[c]{c}%
\quad\quad\quad\quad
\end{array}
\leq\dfrac{1}{\lambda^{1-2\alpha}}\dfrac{4}{1-2\alpha}$

\hfill
\end{proof}

\subsection{Assumptions and notations}

Consider the equation on $\mathbb{R}^{d}$
\begin{equation}
X_{s}^{i}=X_{0}^{i}+\int_{0}^{s}b^{i}(r,X_{r})dr+\sum_{j=1}^{m}\int_{0}%
^{s}\sigma^{i,j}\left(  r,X_{r}\right)  dB_{r}^{j},\,\;s\in\left[  0,T\right]
, \label{se}%
\end{equation}
$i=1,...,d$, where the processes $B^{j}$, $j=1,...,k$ are independent
fractional Brownian motions with Hurst parameter $H$ defined in a complete
probability space $(\Omega,\mathcal{F},\mathbb{P)}$, $X_{0}\ $is a $d$ -
dimensional random variable, and the coefficients $\sigma^{i,j},b^{i}%
:\Omega\times\left[  0,T\right]  \times\mathbb{R}^{d}\mathbb{\rightarrow
}\mathbb{R}$ are measurable functions. Setting $\sigma=\left(  \sigma
^{i,j}\right)  _{d\times k},\,\;\,b=\left(  b^{i}\right)  _{d\times1},$
$B_{s}=\left(  B_{s}^{j}\right)  _{k\times1}$ and $X_{s}=\left(  X_{s}%
^{i}\right)  _{d\times1}$ then we can write the equation (\ref{se}) in a
simpler form%
\[
X_{s}=X_{0}+\int_{0}^{s}b(r,X_{r})dr+\int_{0}^{s}\sigma\left(  r,X_{r}\right)
dB_{r},\,\;s\in\left[  0,T\right]  .
\]

Let us consider the following assumptions on the coefficients, which are
supposed to hold for $\mathbb{P}$ - almost all $\omega\in\Omega$. The
constants $M_{0}$, $M_{R}$, $L_{0}$, $L_{R}$ may depend on $\omega$.

\begin{description}
\item[$\left(  \mathbf{H}_{1}\right)  $] $\sigma(t,x)$ is differentiable in
$x$, and \ there exist some constants $\beta,\delta$, $0<\beta,\delta\leq1,$
and for every $R\geq0$ there exists $M_{R}>0$ such that the following
properties hold for all $t\in\left[  0,T\right]  ,$ $\mathbb{P}-a.s.\;\omega
\in\Omega:$%
\[
\left(  H_{\sigma}\right)  :\left\{
\begin{array}
[c]{rl}%
i)\  & \left\vert \sigma(t,x)-\sigma(s,y)\right\vert \leq M_{0}\left(
|t-s|^{\beta}+|x-y|\right)  ,\quad\forall x,y\in\mathbb{R}^{d},\medskip\\
ii)\  & |\nabla_{x}\sigma(t,y)-\nabla_{x}\sigma(s,z)|\leq M_{R}\left(
|t-s|^{\beta}+|y-z|^{\delta}\right)  ,\ \forall\left\vert y\right\vert
,\left\vert z\right\vert \leq R,
\end{array}
\right.
\]
where $\nabla_{x}\sigma(t,x)=\left(  \nabla_{x}\sigma^{i,j}\left(  t,x\right)
\right)  _{i=\overline{1,d},~j=\overline{1,k}}$ and%
\[
\left\vert \nabla_{x}\sigma(t,x)\right\vert ^{2}=%
%TCIMACRO{\dsum \limits_{\ell=1}^{d}}%
%BeginExpansion
{\displaystyle\sum\limits_{\ell=1}^{d}}
%EndExpansion%
%TCIMACRO{\dsum \limits_{i=1}^{d}}%
%BeginExpansion
{\displaystyle\sum\limits_{i=1}^{d}}
%EndExpansion%
%TCIMACRO{\dsum \limits_{j=1}^{k}}%
%BeginExpansion
{\displaystyle\sum\limits_{j=1}^{k}}
%EndExpansion
\left\vert \partial_{x_{\ell}}\sigma^{i,j}\left(  t,x\right)  \right\vert ^{2}%
\]

\end{description}

\noindent Remark that for all $x\in\mathbb{R}^{d}$
\[
\left\vert \sigma\left(  t,x\right)  \right\vert \leq\left\vert \sigma\left(
0,0\right)  \right\vert +M_{0}\left(  |t|^{\beta}+|x|\right)  \leq
M_{0,T}(1+|x|)
\]
where $M_{0,T}=\left\vert \sigma\left(  0,0\right)  \right\vert +M_{0}%
+M_{0}T.$

Let%
\[
\alpha_{0}=\min\left\{  \frac{1}{2},\beta,\frac{\delta}{1+\delta}\right\}  .
\]
With respect to the coefficient $b$ we assume

\begin{description}
\item[$\left(  \mathbf{H}_{2}\right)  $] There exist $\mu\in(1-\alpha_{0},1]$
and for every $R\geq0$ there exists $L_{R}>0$ such that the following
properties hold for all $t\in\left[  0,T\right]  ,$ $\mathbb{P}-a.s.\;\omega
\in\Omega:$%
\[
\left(  H_{b}\right)  :\left\{
\begin{array}
[c]{rl}%
i)\quad & \left\vert b(r,x)-b(s,y)\right\vert \leq L_{R}\left(  |r-s|^{\mu
}+|x-y|\right)  ,\quad\forall\left\vert x\right\vert ,\left\vert y\right\vert
\leq R,\smallskip\medskip\\
ii)\quad & \left\vert b(t,x)\right\vert \,\leq\,L_{0}(1+|x|),\;\forall
x\in\mathbb{R}^{d}.
\end{array}
\right.
\]

\end{description}

From the work of D. Nualart and\ A. R\u{a}\c{s}canu \cite{nr} we deduce that
under the assumptions $\left(  H_{1}\right)  $ and $\left(  H_{2}\right)  $
with $\beta>1-H$ and $\delta>\frac{1}{H}-1$ and for every fixed $\left(
t,\xi\right)  \in\left[  0,T\right]  \times L^{0}\left(  \Omega,\mathcal{F}%
,\mathbb{P\,};\mathbb{R}^{d})\right)  ,$ the SDE (\ref{se}) has a unique
solution $X^{t,\xi}\in L^{0}\left(  \Omega,\mathcal{F},\mathbb{P\,}%
;W^{\alpha,\infty}(t,T;\mathbb{R}^{d})\right)  ,$ for all $\alpha\in\left(
1-H,\alpha_{0}\right)  .$ Moreover, for $\mathbb{P}$ - almost all $\omega
\in\Omega$
\[
X\left(  \omega,\cdot\right)  =\left(  X^{i}\left(  \omega,\cdot\right)
\right)  _{d\times1}\in C^{1-\alpha}\left(  0,T;\mathbb{R}^{d}\right)  .
\]
In this paper we develop an other proof and we obtain a stronger existence
result: if the starting point $\xi$ belongs to a closed set $K\subset
\mathbb{R}^{d}$ then there exists a unique solution $X^{t,\xi}$ and the
solution evolves in $K.$

\subsection{Basic estimates}

\label{Section Estimates}

Fix a parameter $0<\alpha<\dfrac{1}{2}\wedge\beta$. Given two functions $f\in
W^{\alpha,1}(t,T;\mathbb{R}^{d})$ and $g\in\tilde{W}^{1-\alpha,\infty
}(t,T;\mathbb{R}^{k})$ we denote
\[
G_{t,s}^{\left(  \sigma\right)  }\left(  f\right)  =\int_{t}^{s}\sigma\left(
r,f\left(  r\right)  \right)  dg\left(  r\right)  .
\]
From Proposition \ref{p-est} we infer

\begin{corollary}
\label{c-est-sig}Let the assumption $\left(  H_{1}\right)  $ be satisfied and
$0<\alpha<\frac{1}{2}\wedge\beta$. Let $f\in W^{\alpha,\infty}\left(
t,T;\mathbb{R}^{d}\right)  $ and $g\in\tilde{W}^{1-\alpha,\infty
}(t,T;\mathbb{R}^{k}).$ Then $G_{t,\cdot}^{\left(  \sigma\right)  }\left(
f\right)  \in C^{1-\alpha}\left(  t,T;\mathbb{R}^{d}\right)  $ and for all
$\lambda\geq1$
\begin{equation}%
\begin{array}
[c]{rl}%
i)\  & \left\Vert G_{t,\cdot}^{\left(  \sigma\right)  }\left(  f\right)
\right\Vert _{1-\alpha;\left[  t,T\right]  }\leq C_{0}^{\left(  \sigma
1\right)  }\Lambda_{\alpha}(g;\left[  t,T\right]  )\left(  1+\left\Vert
f\right\Vert _{\alpha,\infty;\left[  t,T\right]  }\right)  ,\medskip\\
ii)\  & \left\vert \!\left\vert \!\left\vert G_{t,\cdot}^{\left(
\sigma\right)  }\left(  f\right)  \right\vert \!\right\vert \!\right\vert
_{\alpha,\lambda;\left[  t,T\right]  }\leq\dfrac{C_{0}^{\left(  \sigma
2\right)  }\Lambda_{\alpha}(g;\left[  t,T\right]  )}{\lambda^{1-2\alpha}%
}\left(  1+\left\Vert f\right\Vert _{\alpha,\lambda;\left[  t,T\right]
}\right)  ,
\end{array}
\label{est-l1-1}%
\end{equation}
where $C_{0}^{\left(  \sigma1\right)  }$ and $C_{0}^{\left(  \sigma2\right)
}$ are constants which only depend on $M_{0,T}$, $M_{0}$, $T$, $\alpha$,
$\beta.$
\end{corollary}

\begin{proof}
From Proposition \ref{p-est} we have%
\[
\left\Vert G_{t,\cdot}^{\left(  \sigma\right)  }\left(  f\right)  \right\Vert
_{1-\alpha;\left[  t,T\right]  }=\left\Vert G_{t,\cdot}\left(  \sigma\left(
\cdot,f\left(  \cdot\right)  \right)  \right)  \right\Vert _{1-\alpha;\left[
t,T\right]  }\leq A_{\alpha,T}^{\left(  1\right)  }~\Lambda_{\alpha}(g;\left[
t,T\right]  )\left\Vert \sigma\left(  \cdot,f\left(  \cdot\right)  \right)
\right\Vert _{\alpha,\infty;\left[  t,T\right]  }%
\]
and%
\[
\left\vert \!\left\vert \!\left\vert G_{t,\cdot}^{\left(  \sigma\right)
}\left(  f\right)  \right\vert \!\right\vert \!\right\vert _{\alpha
,\lambda;\left[  t,T\right]  }\leq\dfrac{\Lambda_{\alpha}(g;\left[
t,T\right]  )}{\lambda^{1-2\alpha}}A_{\alpha,T}^{\left(  2\right)
}~\left\Vert \sigma\left(  \cdot,f\left(  \cdot\right)  \right)  \right\Vert
_{\alpha,\lambda;\left[  t,T\right]  }%
\]
Now the inequalities (\ref{est-l1-1}-i, ii) clearly follow, since for all
$\lambda\geq0,$%
\begin{align*}
&  \left\Vert \sigma\left(  \cdot,f\left(  \cdot\right)  \right)  \right\Vert
_{\alpha,\lambda;\left[  t,T\right]  }\\
&  =\sup_{r\in\left[  t,T\right]  }e^{-\lambda r}\left[  \left\vert
\sigma\left(  r,f\left(  r\right)  \right)  \right\vert +\int_{t}^{r}%
\frac{\left\vert \sigma\left(  r,f\left(  r\right)  \right)  -\sigma\left(
s,f\left(  s\right)  \right)  \right\vert }{\left(  r-s\right)  ^{1+\alpha}%
}ds\right] \\
&  \leq\sup_{r\in\left[  t,T\right]  }e^{-\lambda r}\left[  M_{0,T}\left(
1+\left\vert f\left(  r\right)  \right\vert \right)  +M_{0}\int_{t}^{r}%
\frac{\left(  r-s\right)  ^{\beta}+\left\vert f\left(  r\right)  -f\left(
s\right)  \right\vert }{\left(  r-s\right)  ^{1+\alpha}}ds\right] \\
&  \leq\left(  M_{0,T}+M_{0}\right)  \left(  1+\dfrac{T^{\beta-\alpha}}%
{\beta-\alpha}\right)  \left(  1+\left\Vert f\right\Vert _{\alpha
,\lambda;\left[  t,T\right]  }\right)
\end{align*}

\hfill
\end{proof}

\begin{lemma}
\label{l1}Let the assumption $\left(  H_{1}\right)  $ be satisfied and
$0<\alpha<\frac{1}{2}\wedge\beta$. Let $f,h\in W^{\alpha,\infty}\left(
t,T;\mathbb{R}^{d}\right)  $ and $g\in\tilde{W}^{1-\alpha,\infty
}(t,T;\mathbb{R}^{k}).$ Then for $\left\Vert f\right\Vert _{\infty;\left[
t,T\right]  }\leq R,$ $\left\Vert h\right\Vert _{\infty;\left[  t,T\right]
}\leq R$ it follows%
\begin{equation}%
\begin{array}
[c]{r}%
\left\Vert G_{t,\cdot}^{\left(  \sigma\right)  }\left(  f\right)  -G_{t,\cdot
}^{\left(  \sigma\right)  }\left(  h\right)  \right\Vert _{\alpha
,\lambda;\left[  t,T\right]  }\leq\dfrac{C_{R}^{(\sigma3)}\Lambda_{\alpha
}(g;\left[  t,T\right]  )}{\lambda^{1-2\alpha}}\quad\quad\quad\quad\quad
\quad\medskip\smallskip\\
\times\left(  1+\Delta_{\left[  t,T\right]  }\left(  f\right)  +\Delta
_{\left[  t,T\right]  }\left(  h\right)  \right)  \left\Vert f-h\right\Vert
_{\alpha,\lambda;\left[  t,T\right]  }%
\end{array}
\label{est-l1-2}%
\end{equation}
for all $\lambda\geq1,$ where%
\[
\Delta_{\left[  t,T\right]  }\left(  f\right)  =\sup_{r\in\lbrack t,T]}%
\int_{t}^{r}\frac{|f_{r}-f_{s}|^{\delta}}{(r-s)^{\alpha+1}}\ ds,
\]
and $C_{R}^{(\sigma3)}$ is a constant only depending of $M_{0}$, $M_{R}$, $T$,
$\alpha$, $\beta.$
\end{lemma}

\begin{remark}
\label{delta} If $0<\alpha<\frac{\delta}{1+\delta}$, $\left\Vert f\right\Vert
_{\infty;\left[  t,T\right]  }\leq R$ and $\left\Vert f\right\Vert
_{1-\alpha;\left[  t,T\right]  }\leq C_{R}$~,\ with $C_{R}$ a constant
independent of $t$, then%
\[
\Delta_{\left[  t,T\right]  }\left(  f\right)  \leq C_{R}\sup_{r\in\lbrack
t,T]}\int_{t}^{r}\frac{(r-s)^{\left(  1-\alpha\right)  \delta}}{(r-s)^{\alpha
+1}}ds\leq\dfrac{T^{\delta-\alpha\left(  1+\delta\right)  }}{\delta
-\alpha\left(  1+\delta\right)  }C_{R}~.
\]

\end{remark}

\begin{proof}
[Proof of Lemma \ref{l1}]For the proof we use the ideas from \cite{nr}. By the
inequality (\ref{est11}) we have%
\begin{equation}
\left\Vert G_{t,\cdot}^{\left(  \sigma\right)  }\left(  f\right)  -G_{t,\cdot
}^{\left(  \sigma\right)  }\left(  h\right)  \right\Vert _{\alpha
,\lambda;\left[  t,T\right]  }\leq\frac{\Lambda_{\alpha}(g;\left[  t,T\right]
)A_{\alpha,T}^{\left(  2\right)  }}{\lambda^{1-2\alpha}}\left\Vert
\sigma\left(  \cdot,f\left(  \cdot\right)  \right)  -\sigma\left(
\cdot,h\left(  \cdot\right)  \right)  \right\Vert _{\alpha,\lambda;\left[
t,T\right]  }. \label{est-fh}%
\end{equation}
Remark that if $\left\vert x\right\vert ,\left\vert y\right\vert ,\left\vert
u\right\vert ,\left\vert v\right\vert \leq R$, then%
\[%
\begin{array}
[c]{l}%
\left\vert \sigma^{i,j}\left(  r,x\right)  -\sigma^{i,j}\left(  r,u\right)
-\sigma^{i,j}\left(  s,y\right)  +\sigma^{i,j}\left(  s,v\right)  \right\vert
\smallskip\\
\quad=\left\vert
%TCIMACRO{\dint _{0}^{1}}%
%BeginExpansion
{\displaystyle\int_{0}^{1}}
%EndExpansion
\left\langle x-u,\nabla_{x}\sigma^{i,j}(r,\theta x+(1-\theta)u\right\rangle
d\theta-%
%TCIMACRO{\dint _{0}^{1}}%
%BeginExpansion
{\displaystyle\int_{0}^{1}}
%EndExpansion
\left\langle y-v,\nabla_{x}\sigma^{i,j}(s,\theta y+(1-\theta)v)\right\rangle
d\theta\right\vert \\
\quad\leq\left\vert
%TCIMACRO{\dint _{0}^{1}}%
%BeginExpansion
{\displaystyle\int_{0}^{1}}
%EndExpansion
\left\langle x-y-u+v,\nabla_{x}\sigma^{i,j}(s,\theta y+(1-\theta
)v)\right\rangle d\theta\right\vert \\
\quad\quad\quad+\left\vert
%TCIMACRO{\dint _{0}^{1}}%
%BeginExpansion
{\displaystyle\int_{0}^{1}}
%EndExpansion
\left\langle x-u,\nabla_{x}\sigma^{i,j}(r,\theta x+(1-\theta)u)-\nabla
_{x}\sigma^{i,j}(s,\theta y+(1-\theta)v)\right\rangle d\theta\right\vert
\end{array}
\]
and therefore%
\begin{align*}
\left\vert \sigma\left(  r,x\right)  -\sigma\left(  r,u\right)  -\sigma\left(
s,y\right)  +\sigma\left(  s,v\right)  \right\vert  &  \leq M_{0}\left\vert
x-y-u+v\right\vert \\
&  +M_{R}\left\vert x-u\right\vert \left(  \left\vert s-r\right\vert ^{\beta
}+\left\vert x-y\right\vert ^{\delta}+\left\vert u-v\right\vert ^{\delta
}\right)  .
\end{align*}
Hence for $\left\Vert f\right\Vert _{\infty,\left[  t,T\right]  }\leq R$ and
$\left\Vert h\right\Vert _{\infty,\left[  t,T\right]  }\leq R$\medskip
\newline$\left\Vert \sigma\left(  \cdot,f\left(  \cdot\right)  \right)
-\sigma\left(  \cdot,h\left(  \cdot\right)  \right)  \right\Vert
_{\alpha,\lambda;\left[  t,T\right]  }$\medskip\newline$%
\begin{array}
[c]{c}%
\quad\quad
\end{array}
=\sup\limits_{r\in\left[  t,T\right]  }e^{-\lambda r}\bigg[\left\vert
\sigma\left(  r,f\left(  r\right)  \right)  -\sigma\left(  r,h\left(
r\right)  \right)  \right\vert \medskip\newline%
\begin{array}
[c]{c}%
\quad\quad\quad\quad\quad\quad
\end{array}
+%
%TCIMACRO{\dint _{t}^{r}}%
%BeginExpansion
{\displaystyle\int_{t}^{r}}
%EndExpansion
\dfrac{\left\vert \sigma\left(  r,f\left(  r\right)  \right)  -\sigma\left(
r,h\left(  r\right)  \right)  -\sigma\left(  s,f\left(  s\right)  \right)
+\sigma\left(  s,h\left(  s\right)  \right)  \right\vert }{\left(  r-s\right)
^{1+\alpha}}ds$\medskip\newline$%
\begin{array}
[c]{c}%
\quad\quad
\end{array}
\leq\left(  M_{0}+M_{R}\right)  \sup\limits_{r\in\left[  t,T\right]
}e^{-\lambda r}\bigg[\left\vert f\left(  r\right)  -h\left(  r\right)
\right\vert +%
%TCIMACRO{\dint _{t}^{r}}%
%BeginExpansion
{\displaystyle\int_{t}^{r}}
%EndExpansion
\dfrac{\left\vert f\left(  r\right)  -f\left(  s\right)  -h\left(  r\right)
+h\left(  s\right)  \right\vert }{\left(  r-s\right)  ^{1+\alpha}}ds$%
\medskip\newline$%
\begin{array}
[c]{c}%
\quad\quad\quad\quad\quad\quad
\end{array}
+\left\vert f\left(  r\right)  -h\left(  r\right)  \right\vert \left(
\dfrac{1}{\beta-\alpha}\left(  r-t\right)  ^{\beta-\alpha}+\Delta_{\left[
t,T\right]  }\left(  f\right)  +\Delta_{\left[  t,T\right]  }\left(  h\right)
\right)  \bigg],$\medskip\newline where%
\[
\Delta_{\left[  t,T\right]  }\left(  f\right)  =\sup_{r\in\left[  t,T\right]
}\int_{t}^{r}\frac{\left\vert f\left(  r\right)  -f\left(  s\right)
\right\vert ^{\delta}}{\left(  r-s\right)  ^{1+\alpha}}%
\]
and similar for $h.$ We conclude that%
\begin{align*}
&  \left\Vert \sigma\left(  \cdot,f\left(  \cdot\right)  \right)
-\sigma\left(  \cdot,h\left(  \cdot\right)  \right)  \right\Vert
_{\alpha,\lambda;\left[  t,T\right]  }\\
&  \leq\left(  M_{0}+M_{R}\right)  \left(  1+\frac{T^{\beta-\alpha}}%
{\beta-\alpha}\right)  \left[  1+\Delta_{\left[  t,T\right]  }\left(
f+\Delta_{\left[  t,T\right]  }\left(  h\right)  \right)  \right]  \left\Vert
f-h\right\Vert _{\alpha,\lambda;\left[  t,T\right]  }%
\end{align*}
and the inequality (\ref{est-l1-2}) now follows from (\ref{est-fh}).

\hfill
\end{proof}

\bigskip

We also give similar estimates for%
\[
F_{t,s}^{\left(  b\right)  }\left(  f\right)  =\int_{t}^{s}b\left(
r,f(r)\right)  dr.
\]
where $b$ satisfies the assumptions $\left(  H_{2}\right)  $. Very similar
estimates are given in the paper of Nualart\&R\u{a}\c{s}canu \cite{nr}.

\begin{lemma}
\label{l2}Let $f\in W^{\alpha,\infty}(t,T;\mathbb{R}^{d}).\;$Then $F_{t,\cdot
}^{\left(  b\right)  }\left(  f\right)  =\int_{t}^{\cdot}b\left(
y,f(y)\right)  dy\in C^{1-\alpha}\left(  t,T;\mathbb{R}^{d}\right)  $ and for
all $\lambda\geq1:$
\begin{equation}%
\begin{array}
[c]{rl}%
\left(  j\right)  \; & \left\Vert F_{t,\cdot}^{\left(  b\right)  }\left(
f\right)  \right\Vert _{1-\alpha;\left[  t,T\right]  }\leq C_{0}^{\left(
b1\right)  }\left(  1+\left\Vert f\right\Vert _{\infty;\left[  t,T\right]
}\right)  ,\smallskip\\
\left(  jj\right)  \; & \left\Vert F_{t,\cdot}^{\left(  b\right)  }\left(
f\right)  \right\Vert _{\alpha,\lambda;\left[  t,T\right]  }\leq\dfrac
{C_{0}^{\left(  b2\right)  }}{\lambda^{\alpha}}\left(  1+\left\Vert
f\right\Vert _{\alpha,\lambda;\left[  t,T\right]  }\right)  ,
\end{array}
\label{est-b1}%
\end{equation}
where $C_{0}^{\left(  b1\right)  }$ and $C_{0}^{\left(  b2\right)  }$ are
positive constants depending only on $\alpha,T$ and $L_{0}~.$\newline If
$f,h\in W^{\alpha,\infty}(t,T;\mathbb{R}^{d})$ such that $\ \left\Vert
f\right\Vert _{\infty;\left[  t,T\right]  }\,\leq R$, $\left\Vert h\right\Vert
_{\infty;\left[  t,T\right]  }\leq R$, $\,\,$then%
\begin{equation}
\left\Vert F_{t,\cdot}^{\left(  b\right)  }\left(  f\right)  -F_{t,\cdot
}^{\left(  b\right)  }\left(  h\right)  \right\Vert _{\alpha,\lambda;\left[
t,T\right]  }\leq\dfrac{C_{R}^{(b3)}}{\lambda^{\alpha}}\left\Vert
f-h\right\Vert _{\alpha,\lambda;\left[  t,T\right]  } \label{est-b2}%
\end{equation}
for all $\lambda\geq1$, where $C_{R}^{(b3)}$ constants depending only on
$\alpha,T$ and $L_{R}$ from $\left(  H_{2}\right)  $.
\end{lemma}

\begin{proof}
It is easy to see that $F_{t,\cdot}^{\left(  b\right)  }\left(  f\right)  \in
C^{1}\left(  \left[  t,T\right]  \right)  $ and for $t\leq r\leq s\leq T$%
\[
\left\vert F_{t,s}^{\left(  b\right)  }\left(  f\right)  -F_{t,r}^{\left(
b\right)  }\left(  f\right)  \right\vert =\left\vert F_{r,s}^{\left(
b\right)  }\left(  f\right)  \right\vert \leq L_{0}\left(  1+\left\Vert
f\right\Vert _{\infty;\left[  t,T\right]  }\right)  \left(  s-r\right)  .
\]
Hence the inequality (\ref{est-b1}-j) follows with $C_{0}^{\left(  b1\right)
}=L_{0}\left(  T+T^{\alpha}\right)  .$

Denoting%
\[
F_{t,s}\left(  f\right)  =%
%TCIMACRO{\dint _{t}^{s}}%
%BeginExpansion
{\displaystyle\int_{t}^{s}}
%EndExpansion
f\left(  r\right)  dr
\]
we remark that\medskip\newline$\left\vert F_{t,s}\left(  f\right)  \right\vert
+%
%TCIMACRO{\dint _{t}^{s}}%
%BeginExpansion
{\displaystyle\int_{t}^{s}}
%EndExpansion
\dfrac{|F_{t,s}\left(  f\right)  -F_{t,r}\left(  f\right)  |}{(s-r)^{\alpha
+1}}\ dr\,\smallskip\newline%
\begin{array}
[c]{c}%
\quad\quad\quad
\end{array}
\leq%
%TCIMACRO{\dint _{t}^{s}}%
%BeginExpansion
{\displaystyle\int_{t}^{s}}
%EndExpansion
\left\vert f(r)\right\vert dr+%
%TCIMACRO{\dint _{t}^{s}}%
%BeginExpansion
{\displaystyle\int_{t}^{s}}
%EndExpansion
(s-r)^{-\alpha-1}\left(  \int_{r}^{s}\left\vert f(u)\right\vert du\right)
dr\smallskip$\newline$%
\begin{array}
[c]{c}%
\quad\quad\quad
\end{array}
\leq%
%TCIMACRO{\dint _{t}^{s}}%
%BeginExpansion
{\displaystyle\int_{t}^{s}}
%EndExpansion
\left\vert f(r)\right\vert dr+\dfrac{1}{\alpha}%
%TCIMACRO{\dint _{t}^{s}}%
%BeginExpansion
{\displaystyle\int_{t}^{s}}
%EndExpansion
(s-t)^{-\alpha}\ \ \left\vert f(u)\right\vert du\smallskip$\newline$%
\begin{array}
[c]{c}%
\quad\quad\quad
\end{array}
\leq\left(  T^{\alpha}+\dfrac{1}{\alpha}\right)  \left(  s-t\right)
^{-\alpha}%
%TCIMACRO{\dint _{t}^{s}}%
%BeginExpansion
{\displaystyle\int_{t}^{s}}
%EndExpansion
\left\vert f(r)\right\vert dr,\smallskip$\newline and therefore%
\begin{align*}
\left\vert \!\left\vert \!\left\vert F_{t,\cdot}^{\left(  b\right)  }\left(
f\right)  \right\vert \!\right\vert \!\right\vert _{\alpha,\lambda;\left[
t,T\right]  }  &  \leq\left(  T^{\alpha}+\dfrac{1}{\alpha}\right)  \sup
_{s\in\lbrack t,T]}e^{-\lambda s}\left(  s-t\right)  ^{-\alpha}%
%TCIMACRO{\dint _{t}^{s}}%
%BeginExpansion
{\displaystyle\int_{t}^{s}}
%EndExpansion
\left\vert b(r,f(r))\right\vert dr\\
&  \leq L_{0}\left(  T^{\alpha}+\dfrac{1}{\alpha}\right)  \sup_{s\in\lbrack
t,T]}%
%TCIMACRO{\dint _{t}^{s}}%
%BeginExpansion
{\displaystyle\int_{t}^{s}}
%EndExpansion
e^{-\lambda\left(  s-r\right)  }\frac{1+e^{-\lambda r}\left\vert
f(r))\right\vert }{\left(  s-r\right)  ^{\alpha}}dr,\\
&  \leq L_{0}\left(  T^{\alpha}+\dfrac{1}{\alpha}\right)  \frac{1}%
{\lambda^{\alpha}}\frac{T^{1-2\alpha}}{1-2\alpha}~\left(  1+\left\vert
\!\left\vert \!\left\vert f\right\vert \!\right\vert \!\right\vert
_{\alpha,\lambda;\left[  t,T\right]  }\right)
\end{align*}
that clearly yields (\ref{est-b1}-jj) with $C_{0}^{\left(  b2\right)  }%
=L_{0}\left(  T^{\alpha}+\dfrac{1}{\alpha}\right)  \dfrac{T^{1-2\alpha}%
}{1-2\alpha},$ since $e^{-\lambda\left(  s-r\right)  }\left(  s-r\right)
^{-a}\leq\dfrac{1}{\lambda^{a}}\left(  s-r\right)  ^{-2a}$.

Let now $f,h\in W^{\alpha,\infty}(t,T;\mathbb{R}^{d})$ and $\left\vert
f\right\vert \leq R$ and $\left\vert h\right\vert \leq R.$ Then as here above,
for all $\lambda\geq1,$%
\begin{align*}
\left\Vert F_{t,\cdot}^{\left(  b\right)  }\left(  f\right)  -F_{t,\cdot
}^{\left(  b\right)  }\left(  h\right)  \right\Vert _{\alpha,\lambda;\left[
t,T\right]  }  &  \leq\left(  T^{\alpha}+\dfrac{1}{\alpha}\right)  \sup
_{s\in\lbrack t,T]}e^{-\lambda s}\left(  s-t\right)  ^{-\alpha}\int_{t}%
^{s}\left\vert b(r,f\left(  r\right)  )-b(r,h\left(  r\right)  )\right\vert
dr\\
&  \leq L_{R}\left(  T^{\alpha}+\dfrac{1}{\alpha}\right)  \frac{1}%
{\lambda^{\alpha}}\dfrac{T^{1-2\alpha}}{1-2\alpha}~\left\vert \!\left\vert
\!\left\vert f-h\right\vert \!\right\vert \!\right\vert _{\alpha
,\lambda;\left[  t,T\right]  }~.
\end{align*}
Hence the inequality (\ref{est-b2})holds with $C_{R}^{(b3)}=L_{R}\left(
T^{\alpha}+\dfrac{1}{\alpha}\right)  \dfrac{T^{1-2\alpha}}{1-2\alpha}.$

\hfill
\end{proof}

Finally we present some auxiliary estimates used in the sequel.

\begin{lemma}
\label{aux}Let the assumptions $(H_{1})$ and $(H_{2})$ be satisfied and
$0<\alpha<\beta\wedge\frac{1}{2}$. If $Y$ is a Holder continuous function with
$\left\Vert Y\right\Vert _{1-\alpha;\left[  t,T\right]  }\leq R$ then there
exist some positive constants $C_{R}^{\left(  1\right)  }=\left(
R+1+T\right)  L_{R}$, $C_{R}^{\left(  2\right)  }=C^{\left(  2\right)
}\left(  R,M_{0},T,\alpha,\beta,\Lambda_{\alpha}\left(  g\right)  \right)  $,
$C_{R}^{\left(  3\right)  }=2\left(  1+R\right)  L_{0}$ and $C_{R}^{\left(
4\right)  }=C^{\left(  4\right)  }\left(  R,M_{0},T,\alpha,\beta
,\Lambda_{\alpha}\left(  g\right)  \right)  $ such that for all $0\leq t\leq
s\leq T:$
\begin{equation}%
\begin{array}
[c]{ll}%
\left(  a\right)  \; & \left\vert
%TCIMACRO{\dint _{t}^{s}}%
%BeginExpansion
{\displaystyle\int_{t}^{s}}
%EndExpansion
\left[  b(r,Y_{r})-b\left(  t,Y_{t}\right)  \right]  dr\right\vert \leq
C_{R}^{\left(  1\right)  }\left(  s-t\right)  ^{2-\alpha}\quad\text{and}%
\medskip\\
\left(  b\right)  \; & \left\vert
%TCIMACRO{\dint _{t}^{s}}%
%BeginExpansion
{\displaystyle\int_{t}^{s}}
%EndExpansion
\left[  \sigma(r,Y_{r})-\sigma\left(  t,Y_{t}\right)  \right]  dg\left(
r\right)  \right\vert \leq C_{R}^{\left(  2\right)  }\left(  s-t\right)
^{1+\min\left\{  \beta-\alpha,1-2\alpha\right\}  }.
\end{array}
\label{aux1}%
\end{equation}
and for all $0\leq t\leq\tau\leq s\leq T:$%
\begin{equation}%
\begin{array}
[c]{ll}%
\left(  c\right)  \; & \left\vert
%TCIMACRO{\dint _{\tau}^{s}}%
%BeginExpansion
{\displaystyle\int_{\tau}^{s}}
%EndExpansion
\left[  b(r,Y_{r})-b\left(  t,Y_{t}\right)  \right]  dr\right\vert \leq
C_{R}^{\left(  3\right)  }\left(  s-\tau\right)  \quad\text{and}\medskip\\
\left(  d\right)  \; & \left\vert
%TCIMACRO{\dint _{\tau}^{s}}%
%BeginExpansion
{\displaystyle\int_{\tau}^{s}}
%EndExpansion
\left[  \sigma(r,Y_{r})-\sigma\left(  t,Y_{t}\right)  \right]  dg\left(
r\right)  \right\vert \leq C_{R}^{\left(  4\right)  }\left(  s-\tau\right)
^{1-\alpha}.
\end{array}
\label{aux2}%
\end{equation}

\end{lemma}

\begin{proof}
$\left(  a\right)  \quad$We have%
\begin{align*}
\left\vert \int_{t}^{s}\left[  b(r,Y_{r})-b\left(  t,Y_{t}\right)  \right]
dr\right\vert  &  \leq\left(  s-t\right)  \sup_{r\in\left[  t,s\right]
}\left\vert b(r,Y_{r})-b\left(  t,Y_{t}\right)  \right\vert \\
&  \leq\left(  s-t\right)  \sup_{r\in\left[  t,s\right]  }L_{R}\left(
\left\vert Y_{r}-Y_{t}\right\vert +\left\vert r-t\right\vert ^{\mu}\right) \\
&  \leq\left(  s-t\right)  L_{R}\sup_{r\in\left[  t,s\right]  }\left\vert
\left[  R+\left(  1+T\right)  \right]  \left(  r-t\right)  ^{1-\alpha
}\right\vert \\
&  \leq C_{R}^{\left(  1\right)  }\left(  s-t\right)  ^{2-\alpha}%
\end{align*}
$\left(  b\right)  \quad$By the assumptions $(H_{1})$ and Lemma \ref{l66} we
have
\begin{align*}
\left\vert \sigma(r,Y_{r})-\sigma\left(  \theta,Y_{\theta}\right)
\right\vert  &  \leq M_{0}\left[  \left\vert r-\theta\right\vert ^{\beta
}+\left\vert Y_{r}-Y_{\theta}\right\vert \right] \\
&  \leq M_{0}\left[  \left\vert r-\theta\right\vert ^{\beta}+R\left\vert
r-\theta\right\vert ^{1-\alpha}\right]
\end{align*}
and\medskip\newline$\left\vert
%TCIMACRO{\dint _{t}^{s}}%
%BeginExpansion
{\displaystyle\int_{t}^{s}}
%EndExpansion
\left[  \sigma(r,Y_{r})-\sigma\left(  t,Y_{t}\right)  \right]  dg\left(
r\right)  \right\vert $\medskip\newline$%
\begin{array}
[c]{c}%
\quad\quad\quad
\end{array}
\leq\left\vert \Lambda_{\alpha}\left(  g\right)  \left(
%TCIMACRO{\dint _{t}^{s}}%
%BeginExpansion
{\displaystyle\int_{t}^{s}}
%EndExpansion
\dfrac{\sigma(r,Y_{r})-\sigma\left(  t,Y_{t}\right)  }{\left(  r-t\right)
^{\alpha}}dr+%
%TCIMACRO{\dint _{t}^{s}}%
%BeginExpansion
{\displaystyle\int_{t}^{s}}
%EndExpansion%
%TCIMACRO{\dint _{t}^{r}}%
%BeginExpansion
{\displaystyle\int_{t}^{r}}
%EndExpansion
\dfrac{\sigma(r,Y_{r})-\sigma\left(  \theta,Y_{\theta}\right)  }{\left(
r-\theta\right)  ^{1+\alpha}}d\theta dr\right)  \right\vert $\medskip\newline$%
\begin{array}
[c]{c}%
\quad\quad\quad
\end{array}
\leq\Lambda_{\alpha}\left(  g\right)  M_{0}\left\vert
%TCIMACRO{\dint _{t}^{s}}%
%BeginExpansion
{\displaystyle\int_{t}^{s}}
%EndExpansion
\left[  \left(  r-t\right)  ^{\beta-\alpha}+R\left(  r-t\right)  ^{1-2\alpha
}\right]  dr\right.  \medskip\newline%
\begin{array}
[c]{c}%
\quad\quad\quad\quad\quad\quad\quad\quad
\end{array}
\left.  +%
%TCIMACRO{\dint _{t}^{s}}%
%BeginExpansion
{\displaystyle\int_{t}^{s}}
%EndExpansion%
%TCIMACRO{\dint _{t}^{r}}%
%BeginExpansion
{\displaystyle\int_{t}^{r}}
%EndExpansion
\left[  \left(  r-\theta\right)  ^{\beta-1-\alpha}+R\left(  r-\theta\right)
^{-2\alpha}\right]  d\theta dr\right\vert \medskip$\newline$%
\begin{array}
[c]{c}%
\quad\quad\quad
\end{array}
\leq\tilde{C}_{R}^{\left(  2\right)  }\left\{  \left(  s-t\right)
^{\beta-\alpha+1}+\left(  s-t\right)  ^{2-2\alpha}+\left(  s-t\right)
^{\beta-\alpha+1}+\left(  s-t\right)  ^{2-2\alpha}\right\}  $\medskip\newline$%
\begin{array}
[c]{c}%
\quad\quad\quad
\end{array}
\leq C_{R}^{\left(  2\right)  }\left(  s-t\right)  ^{1+\min\left\{
\beta-\alpha,1-2\alpha\right\}  }.$\medskip\newline$\left(  c\right)  \quad$We
have, by $\left(  H_{2}-ii\right)  ,$%
\begin{align*}
\left\vert
%TCIMACRO{\dint _{\tau}^{s}}%
%BeginExpansion
{\displaystyle\int_{\tau}^{s}}
%EndExpansion
\left[  b(r,Y_{r})-b\left(  t,Y_{t}\right)  \right]  dr\right\vert  &
\leq\left(  s-\tau\right)  \underset{r\in\left[  \tau,s\right]  }{\sup
}\left\vert b(r,Y_{r})-b\left(  t,Y_{t}\right)  \right\vert \\
&  \leq L_{0}\left(  s-\tau\right)  (2+\left\vert Y_{r}\right\vert +\left\vert
Y_{t}\right\vert )\\
&  \leq2\left(  1+R\right)  L_{0}\left(  s-\tau\right)
\end{align*}
$\left(  d\right)  \quad$Using (\ref{est1}) we deduce \newline$\left\vert
%TCIMACRO{\dint _{\tau}^{s}}%
%BeginExpansion
{\displaystyle\int_{\tau}^{s}}
%EndExpansion
\left[  \sigma(r,Y_{r})-\sigma\left(  t,Y_{t}\right)  \right]  dg\left(
r\right)  \right\vert \medskip\newline%
\begin{array}
[c]{c}%
\quad\quad\quad
\end{array}
\leq\Lambda_{\alpha}\left(  g\right)  \left\vert
%TCIMACRO{\dint _{\tau}^{s}}%
%BeginExpansion
{\displaystyle\int_{\tau}^{s}}
%EndExpansion
\dfrac{\sigma(r,Y_{r})-\sigma\left(  t,Y_{t}\right)  }{\left(  r-\tau\right)
^{\alpha}}dr+%
%TCIMACRO{\dint _{\tau}^{s}}%
%BeginExpansion
{\displaystyle\int_{\tau}^{s}}
%EndExpansion%
%TCIMACRO{\dint _{\tau}^{r}}%
%BeginExpansion
{\displaystyle\int_{\tau}^{r}}
%EndExpansion
\dfrac{\sigma(r,Y_{r})-\sigma\left(  \theta,Y_{\theta}\right)  }{\left(
r-\theta\right)  ^{1+\alpha}}d\theta dr\right\vert \medskip\newline%
\begin{array}
[c]{c}%
\quad\quad\quad
\end{array}
\leq\Lambda_{\alpha}\left(  g\right)  M_{0}%
%TCIMACRO{\dint _{\tau}^{s}}%
%BeginExpansion
{\displaystyle\int_{\tau}^{s}}
%EndExpansion
\left(  \dfrac{\left(  r-t\right)  ^{\beta}}{\left(  r-\tau\right)  ^{\alpha}%
}+R\dfrac{\left(  r-t\right)  ^{1-\alpha}}{\left(  r-\tau\right)  ^{\alpha}%
}\right)  dr\medskip\newline%
\begin{array}
[c]{c}%
\quad\quad\quad\quad\quad\quad\quad\quad
\end{array}
+\Lambda_{\alpha}\left(  g\right)  M_{0}%
%TCIMACRO{\dint _{\tau}^{s}}%
%BeginExpansion
{\displaystyle\int_{\tau}^{s}}
%EndExpansion%
%TCIMACRO{\dint _{\tau}^{r}}%
%BeginExpansion
{\displaystyle\int_{\tau}^{r}}
%EndExpansion
\left(  \left(  r-\theta\right)  ^{\beta-\alpha-1}+R\left(  r-\theta\right)
^{-2\alpha}\right)  d\theta dr\medskip\newline%
\begin{array}
[c]{c}%
\quad\quad\quad
\end{array}
\leq\tilde{C}_{R}^{\left(  4\right)  }\left[  T^{\beta}\left(  s-\tau\right)
^{1-\alpha}+T^{1-\alpha}\left(  s-\tau\right)  ^{1-\alpha}+\left(
s-\tau\right)  ^{\beta+1-\alpha}+\left(  s-\tau\right)  ^{2-2\alpha}\right]
\medskip\newline%
\begin{array}
[c]{c}%
\quad\quad\quad
\end{array}
\leq C_{R}^{\left(  4\right)  }\left(  s-\tau\right)  ^{1-\alpha}.$

\hfill
\end{proof}

\section{Fractional viability. Main result.}

Consider the stochastic differential equation \ driven by the fractional
Brownian motion $B$ with Hurst parameter $1/2<H\leq1$, $\mathbb{P}%
-a.s.\;\omega\in\Omega,$%
\begin{equation}
X_{s}^{t,x}=x+\int_{t}^{s\vee t}b(r,X_{r}^{t,x})dr+\int_{t}^{s\vee t}%
\sigma(r,X_{r}^{t,x})dB_{r}^{H},\;\;s\in\lbrack t,T], \label{frac}%
\end{equation}
where

\begin{itemize}
\item $B=\left(  B^{i}\right)  _{k\times1}$ , $B^{i}$, $i=\overline{1,k},$ are
independent fractional Brownian motions with Hurst parameter $H,$ $\dfrac
{1}{2}<H<1$, and the integral with respect to $B$ is a pathwise
Riemann-Stieltjes integral;

\item $X_{0}=\left(  X_{0}^{i}\right)  _{d\times1}$ is a $d$ - dimensional
random variable defined in a complete probability space $(\Omega
,\mathcal{F},\mathbb{P)};$

\item $b:\left[  0,T\right]  \times\mathbb{R}^{d}\mathbb{\rightarrow
}\mathbb{R}^{d}$ and $\sigma:\left[  0,T\right]  \times\mathbb{R}%
^{d}\mathbb{\rightarrow R}^{d\times k}$ are continuous functions.
\end{itemize}

\begin{definition}
Let $\mathcal{K}=\left\{  K\left(  t\right)  :t\in\left[  0,T\right]
\right\}  $ be a family of subsets of $\mathbb{R}^{d}$. We shall say that
$\mathcal{K}$ is viable for the equation (\ref{frac}) if, starting at any time
$t$ $\in\lbrack0,T]$ and from any point $x\in K\left(  t\right)  ,$ at least
one its solution $X_{s}^{t,x}\in K\left(  s\right)  $ for all $s\in\left[
t,T\right]  .$
\end{definition}

\begin{definition}
The family $\mathcal{K}$ is said to be invariant for the equation (\ref{frac})
if, for any $t\in\left[  0,T\right]  $ and for any starting point $x\in$
$K\left(  t\right)  $, all solutions $\left\{  X_{s}^{t,x}:s\in\left[
t,T\right]  \right\}  $ of the fractional stochastic differential equation
(\ref{frac}) have the property%
\[
X_{s}^{t,x}\in K\left(  s\right)  \text{ for all }s\in\left[  t,T\right]  .
\]

\end{definition}

Remark that, in the case when the equation has a unique solution (which is the
case for the equation (\ref{frac}) under the assumptions $(H_{1})$ and
$(H_{2})$), viability is equivalent with invariance.

Assume that the mappings $b$ and $\sigma$ from the equations (\ref{frac}) are
satisfying $(H_{1})$ and $(H_{2})$

\begin{definition}
Let $t\in\left[  0,T\right]  $ and $x\in K\left(  t\right)  .$ Let $\frac
{1}{2}<1-\alpha<H.$ We say that the pair $(b\left(  t,x\right)  ,\sigma\left(
t,x\right)  )$ is $\left(  1-\alpha\right)  -$fractional $B^{H}-$%
\ contingent\ to $K\left(  t\right)  $ in $\left(  t,x\right)  $ if there
exist random variable $\bar{h}=\bar{h}^{t,x}>0,$ a stochastic process
$Q=Q^{t,x}:\Omega\times\left[  t,t+\bar{h}\right]  \rightarrow\mathbb{R}^{d}$
and for every $R>0$ such that $\left\vert x\right\vert \leq R$ there exist two
random variables $H_{R}$,$\tilde{H}_{R}>0$ and a constant $\gamma=\gamma
_{R}\in\left(  0,1\right)  $ which are independent of $\left(  t,\bar
{h}\right)  $ (the constants $H_{R}$, $\tilde{H}_{R}$, $\gamma_{R}$ depend
only on $R$, $L_{R}$, $M_{0,T}$, $M_{0}$, $L_{0}$, $T$, $\alpha$, $\beta$,
$\Lambda_{\alpha}\left(  B^{H}\right)  $) such that for all $s,\tau\in\left[
t,t+\bar{h}\right]  $%
\[
\left\vert Q\left(  s\right)  -Q\left(  \tau\right)  \right\vert \leq
H_{R}\left\vert s-\tau\right\vert ^{1-\alpha}\quad\text{and}\quad\left\vert
Q\left(  s\right)  \right\vert \leq\tilde{H}_{R}\left\vert s-t\right\vert
^{1+\gamma}%
\]
satisfying%
\[
x+\left(  s-t\right)  b\left(  t,x\right)  +\sigma\left(  t,x\right)  \left[
B_{s}^{H}-B_{t}^{H}\right]  +Q\left(  s\right)  \in K\left(  s\right)  .
\]
\newline\hfill
\end{definition}

The main result (existence result and characterization of the viability) of
our paper is the following

\begin{theorem}
\label{main}Let $\mathcal{K}=\left\{  K\left(  t\right)  :t\in\left[
0,T\right]  \right\}  $ be a family of nonempty closed subsets of
$\mathbb{R}^{d}$. Assume that the maps $b$ and $\sigma$ from the equations
(\ref{frac}) are satisfying $(H_{1})$, $(H_{2})$ with $\frac{1}{2}<H<1$,
$1-H<\beta$, $\delta>\frac{1-H}{H}.$ \newline Let $1-H<\alpha<\alpha_{0}$.
Then the following assertions are equivalent:

\begin{itemize}
\item $\mathcal{K}$ is viable for the fractional stochastic differential
equation (\ref{frac}), i.e. for all $t\in\left[  0,T\right]  $ and for all
$x\in K\left(  t\right)  $ there exists a solution $X^{t,x}\left(
\omega,\cdot\right)  \in C^{1-\alpha}\left(  \left[  t,T\right]
;~\mathbb{R}^{d}\right)  $ of the equation
\[
X_{s}^{t,x}=x+\int_{t}^{s\vee t}b(r,X_{r}^{t,x})dr+\int_{t}^{s\vee t}%
\sigma(r,X_{r}^{t,x})dB_{r}^{H},\;\;s\in\lbrack t,T],\ a.s.\;\omega\in\Omega,
\]

and $X_{s}^{t,x}\in K\left(  s\right)  ,$ for all $s\in\left[  t,T\right]  .$

\item For all $t\in\left[  0,T\right]  $ and all $x\in K\left(  t\right)  ,$
$\left(  b\left(  t,x\right)  ,\sigma\left(  t,x\right)  \right)  $ is
$\left(  1-\alpha\right)  -$fractional\ $B^{H}-$contingent\ to $K\left(
t\right)  $ in $\left(  t,x\right)  $ .
\end{itemize}
\end{theorem}

\section{Deterministic approach}

Let arbitrary fixed $\left(  t,x\right)  \in\left[  0,T\right]  \times%
%TCIMACRO{\U{211d} }%
%BeginExpansion
\mathbb{R}
%EndExpansion
^{d}.$ Consider the deterministic differential equation on $%
%TCIMACRO{\U{211d} }%
%BeginExpansion
\mathbb{R}
%EndExpansion
^{d}:$%
\begin{equation}
X_{s}^{tx}=x+\int_{t}^{s}b(r,X_{r}^{tx})dr+\int_{t}^{s}\sigma\left(
r,X_{r}^{tx}\right)  dg\left(  r\right)  ,\,\;s\in\left[  t,T\right]  ,
\label{det}%
\end{equation}
where $g\in\tilde{W}^{1-\alpha,\infty}(t,T;\mathbb{R}^{k})$ and the
coefficients $b:\left[  0,T\right]  \times%
%TCIMACRO{\U{211d} }%
%BeginExpansion
\mathbb{R}
%EndExpansion
^{d}\rightarrow%
%TCIMACRO{\U{211d} }%
%BeginExpansion
\mathbb{R}
%EndExpansion
^{d}$ and $\sigma:\left[  0,T\right]  \times%
%TCIMACRO{\U{211d} }%
%BeginExpansion
\mathbb{R}
%EndExpansion
^{d}\rightarrow\mathbb{R}^{d\times k}$ are continuous functions satisfying the
assumptions $\left(  H_{1}\right)  $, $\left(  H_{2}\right)  $. Let $\alpha$
be arbitrary fixed such that
\[
0<\alpha<\alpha_{0}=\left\{  \frac{1}{2},\beta,\frac{\delta}{1+\delta
}\right\}  .
\]

Nualart and Rascanu proved in \cite{nr} that if the assumptions $\left(
H_{1}\right)  $ and $\left(  H_{2}\right)  $ are satisfied then the equation
(\ref{det}) has a unique solution which is $\left(  1-\alpha\right)  -$ Holder
continuous. In the following Lemma we shall proof that the Holder constant of
this solution has the form $C_{0}\left(  1+\left\vert x\right\vert \right)  ,$
with $C_{0}$ a positive constant depending only on $M_{0,T}$, $M_{0}$, $L_{0}$
$T$, $\alpha$, $\beta.$

\begin{lemma}
\label{l66}Let the assumptions $(H_{1})$ and $(H_{2})$ be satisfied. If
$X^{t,x}$ is a solution of the equation (\ref{det}) then $X^{t,x}$ is $\left(
1-\alpha\right)  $-Holder continuous and
\[
\left\Vert X_{\cdot}^{t,x}\right\Vert _{1-\alpha;\left[  t,T\right]  }\leq
C_{0}\left(  1+\left\vert x\right\vert \right)
\]
where $C_{0}$ is a constant depending only on $M_{0,T}$, $M_{0}$, $L_{0}$ $T$,
$\alpha$, $\beta,$ $\Lambda_{\alpha}\left(  g\right)  .$
\end{lemma}

\begin{proof}
By Corollary \ref{c-est-sig} and Lemma \ref{l2} we have for all $\lambda\geq1$%
\begin{align*}
&  \left\Vert X_{\cdot}^{t,x}\right\Vert _{\alpha,\lambda;\left[  t,T\right]
}\\
&  \leq\left\vert x\right\vert +\left\Vert F_{t,\cdot}^{\left(  b\right)
}\left(  X^{t,x}\right)  \right\Vert _{\alpha,\lambda;\left[  t,T\right]
}+\left\Vert G_{t,\cdot}^{\left(  \sigma\right)  }\left(  X^{t,x}\right)
\right\Vert _{\alpha,\lambda;\left[  t,T\right]  }\medskip\\
&  \leq\left\vert x\right\vert +\dfrac{C_{0}^{\left(  b2\right)  }}%
{\lambda^{1-2\alpha}}\left(  1+\left\Vert X_{\cdot}^{t,x}\right\Vert
_{\alpha,\lambda;\left[  t,T\right]  }\right)  +\frac{\Lambda_{a}\left(
g;\left[  0,T\right]  \right)  C_{0}^{\left(  \sigma2\right)  }}%
{\lambda^{1-2\alpha}}\left(  1+\left\Vert X_{\cdot}^{t,x}\right\Vert
_{\alpha,\lambda;\left[  t,T\right]  }\right)  \medskip\\
&  \leq\left\vert x\right\vert +\frac{1}{2}\left(  1+\left\Vert X^{t,x}%
\right\Vert _{\alpha,\lambda;\left[  t,T\right]  }\right)
\end{align*}
for $\lambda=\lambda_{0}\geq1$ sufficiently large,%
\[
\dfrac{C_{0}^{\left(  b2\right)  }+\Lambda_{a}\left(  g\right)  C_{0}^{\left(
\sigma2\right)  }}{\lambda_{0}^{1-2\alpha}}\leq\frac{1}{2},
\]
(remark that the constant $C_{0}^{\left(  b2\right)  }+\Lambda_{\alpha}\left(
g\right)  C_{0}^{\left(  \sigma2\right)  }$ is independent of $\lambda).$ Then
we have%
\[
\left\Vert X^{t,x}\right\Vert _{\alpha,\lambda;\left[  t,T\right]  }%
\leq2\left(  1+\left\vert x\right\vert \right)  \quad\text{and}\quad\left\Vert
X_{\cdot}^{t,x}\right\Vert _{\alpha,\infty;\left[  t,T\right]  }\leq2\left(
1+\left\vert x\right\vert \right)  e^{\lambda_{0}T}.
\]
On the other hand, using the same lemmas we get%
\begin{align*}
\left\Vert X^{t,x}\right\Vert _{1-\alpha;\left[  t,T\right]  } &
\leq\left\vert x\right\vert +\left\Vert F_{t,\cdot}^{\left(  b\right)
}\left(  X^{t,x}\right)  \right\Vert _{1-\alpha;\left[  t,T\right]
}+\left\Vert G_{t,\cdot}^{\left(  \sigma\right)  }\left(  X^{t,x}\right)
\right\Vert _{1-\alpha;\left[  t,T\right]  }\smallskip\\
&  \leq\left\vert x\right\vert +C_{0}^{\left(  b1\right)  }\left(
1+\left\Vert X^{t,x}\right\Vert _{\infty;\left[  t,T\right]  }\right)
+\Lambda_{\alpha}\left(  g\right)  C_{0}^{\left(  \sigma1\right)  }\left(
1+\left\Vert X^{t,x}\right\Vert _{\alpha,\infty;\left[  t,T\right]  }\right)
\smallskip\\
&  \leq\left\vert x\right\vert +\left(  C_{0}^{\left(  b1\right)  }%
+\Lambda_{\alpha}\left(  g\right)  C_{0}^{\left(  \sigma1\right)  }\right)
\left(  1+\left\Vert X^{t,x}\right\Vert _{\alpha,\infty;\left[  t,T\right]
}\right)  \smallskip\\
&  \leq\left\vert x\right\vert +\left(  C_{0}^{\left(  b1\right)  }%
+\Lambda_{\alpha}\left(  g\right)  C_{0}^{\left(  \sigma1\right)  }\right)
\left(  1+2e^{\lambda_{0}T}+2\left\vert x\right\vert e^{\lambda_{0}T}\right)
\\
&  \leq C_{0}\left(  1+\left\vert x\right\vert \right)
\end{align*}
Hence $X_{\cdot}^{t,x}$ is $\left(  1-\alpha\right)  -$Holder continuous with
the Holder constant $C_{0}\left(  1+\left\vert x\right\vert \right)  $ where
$C_{0}$ depends only on $M_{0,T}$, $M_{0}$, $L_{0}$ $T$, $\alpha$, $\beta$ and
$\Lambda_{\alpha}\left(  g\right)  .$ The proof is now complete.

\hfill
\end{proof}

\medskip

Assume that the maps $b$ and $\sigma$ from the equations (\ref{det}) are
satisfying $(H_{1})$ and $(H_{2})$.

\begin{definition}
[Tangency property]\label{def-g-tg}Let $t\in\left[  0,T\right]  $ and $x\in
K\left(  t\right)  .$  We say that the pair $(b\left(  t,x\right)
,\sigma\left(  t,x\right)  )$ is $\left(  1-\alpha\right)  -$fractional
$g-$\ tangent\ to $K\left(  t\right)  $ in $\left(  t,x\right)  $ if there
exist $\bar{h}=\bar{h}^{t,x}>0,$ and two functions $U=U^{t,x}:\left[
t,t+\bar{h}\right]  \rightarrow\mathbb{R}^{d}$, $U\left(  t\right)  =0,$ and
$V=V^{t,x}:\left[  t,t+\bar{h}\right]  \rightarrow\mathbb{R}^{d\times k}$,
$V\left(  t\right)  =0,$ and for every $R>0$ such that $\left\vert
x\right\vert \leq R$ there exist two constants $D_{R},\tilde{D}_{R}>0$
independent of $\left(  t,\bar{h}\right)  $ such that for all $s,\tau
\in\left[  t,t+\bar{h}\right]  $%
\begin{equation}
\left\vert U\left(  \tau\right)  -U\left(  s\right)  \right\vert \leq
D_{R}\left\vert \tau-s\right\vert ^{1-\alpha}\quad\text{and}\quad\left\vert
V\left(  \tau\right)  -V\left(  s\right)  \right\vert \leq\tilde{D}%
_{R}\left\vert \tau-s\right\vert ^{\min\left\{  \beta,~1-\alpha\right\}
}\label{UV}%
\end{equation}
and satisfying%
\[
x+\int_{t}^{s}\left[  b(t,x)+U\left(  r\right)  \right]  dr+\int_{t}%
^{s}\left[  \sigma(t,x)+V\left(  r\right)  \right]  dg\left(  r\right)  \in
K\left(  s\right)  ,~\ \ ~\text{for all }s\in\left[  t,t+\bar{h}\right]  .
\]

\end{definition}

\begin{definition}
[Contingency property]\label{def-g-ctg}Let $t\in\left[  0,T\right]  $ and
$x\in K\left(  t\right)  .$ We say that the pair $(b\left(  t,x\right)
,\sigma\left(  t,x\right)  )$ is $\left(  1-\alpha\right)  -$fractional
$g-$\ contingent\ to $K\left(  t\right)  $ in $\left(  t,x\right)  $ if there
exist $\bar{h}=\bar{h}^{t,x}>0,$ a function $Q=Q^{t,x}:\left[  t,t+\bar
{h}\right]  \rightarrow\mathbb{R}^{d}$ and for every $R>0$ such that
$\left\vert x\right\vert \leq R$ there exist two constants $G_{R}$,$\tilde
{G}_{R}>0$ independent of $\left(  t,\bar{h}\right)  $ and a constant
$\gamma=\gamma_{R}\in\left(  0,1\right)  $ also independent of $\left(
t,\bar{h}\right)  $ (the constants $G_{R}$,$\tilde{G}_{R}$, $\gamma_{R}$
depend only on $R,L_{R},M_{0,T}$, $M_{0}$, $L_{0}$ $T$, $\alpha$, $\beta$ and
$\Lambda_{\alpha}\left(  g\right)  $) such that for all $s,\tau\in\left[
t,t+\bar{h}\right]  $%
\[
\left\vert Q\left(  \tau\right)  -Q\left(  s\right)  \right\vert \leq
G_{R}\left\vert \tau-s\right\vert ^{1-\alpha}\quad\text{and}\quad\left\vert
Q\left(  s\right)  \right\vert \leq\tilde{G}_{R}\left\vert s-t\right\vert
^{1+\gamma}%
\]
and satisfying%
\[
x+\left(  s-t\right)  b\left(  t,x\right)  +\sigma\left(  t,x\right)  \left[
g\left(  s\right)  -g\left(  t\right)  \right]  +Q\left(  s\right)  \in
K\left(  s\right)  ,~\ \ ~\text{for all }s\in\left[  t,t+\bar{h}\right]  .
\]

\end{definition}

We can now state the main result of the section.

\begin{theorem}
\label{detmain}Let $\mathcal{K}=\left\{  K\left(  t\right)  :t\in\left[
0,T\right]  \right\}  $ be a family of nonempty closed subsets of
$\mathbb{R}^{d}$. Assume $(H_{1})$ and $(H_{2})$ are satisfied and%
\[
1-\mu<\alpha<\alpha_{0}=\min\left\{  \frac{1}{2},\beta,\frac{\delta}{1+\delta
}\right\}  .
\]
Then the following assertions are equivalent:

\begin{itemize}
\item[$\left(  j\right)  $] $\ \mathcal{K}$ is $C^{1-\alpha}-$viable for the
fractional differential equation (\ref{det}), i.e. for any $t\in\left[
0,T\right]  $ and for any starting point $x\in$ $K\left(  t\right)  $, there
exists a solution $X^{t,x}\left(  \cdot\right)  \in C^{1-\alpha}\left(
\left[  t,T\right]  ;~\mathbb{R}^{d}\right)  $ of the equation%
\begin{equation}
X_{s}^{t,x}=x+\int_{t}^{s\vee t}b(r,X_{r}^{t,x})dr+\int_{t}^{s\vee t}%
\sigma(r,X_{r}^{t,x})dg\left(  r\right)  ,\ s\in\lbrack t,T],\label{det1}%
\end{equation}
such that $X_{s}^{t,x}\in K\left(  s\right)  ,$ for all $s\in\lbrack t,T].$

\item[$\left(  jj\right)  $] For all $t\in\left[  0,T\right]  $ and all $x\in
K\left(  t\right)  $ the pair $\left(  b\left(  t,x\right)  ,\sigma\left(
t,x\right)  \right)  $ is $\left(  1-\alpha\right)  -$fractional
$g-$\ tangent\ to $K\left(  t\right)  $ in $\left(  t,x\right)  .$

\item[$\left(  jjj\right)  $] \ For all $t\in\left[  0,T\right]  $ and all
$x\in K\left(  t\right)  $ the pair $\left(  b\left(  t,x\right)
,\sigma\left(  t,x\right)  \right)  $ is $\left(  1-\alpha\right)
-$fractional $g-$\ contingent\ to $K\left(  t\right)  $ in $\left(
t,x\right)  .$
\end{itemize}
\end{theorem}

\begin{proof}
Let $0<\varepsilon\leq1.$ We denote by $C_{R},$ $C_{R}^{\left(  1\right)  },$
$C_{R}^{\left(  2\right)  }$, $\ldots$ denote a generic positive constant
independent of $\varepsilon$ and depending only on $R,L_{R},M_{0,T}$, $M_{0}$,
$L_{0}$ $T$, $\alpha$, $\beta$ and $\Lambda_{\alpha}\left(  g\right)  .$
\medskip\newline$\left(  \mathbf{j}\right)  \Rightarrow\left(  \mathbf{jj}%
\right)  :$

Let $t\in\left[  0,T\right]  $ and $x\in K\left(  t\right)  $ be arbitrary
fixed and $X^{t,x}\in C^{1-\alpha}\left(  \left[  t,T\right]  ;~\mathbb{R}%
^{d}\right)  $ a solution of the equation (\ref{det1}) such that $X_{s}%
^{t,x}\in K\left(  s\right)  ,$ for all $s\in\lbrack t,T].$ Let $R_{0}>0$ such
that $\left\vert x\right\vert \leq R_{0}.$ Then by Lemma \ref{l66}%
\[
\left\Vert X^{t,x}\right\Vert _{1-\alpha;\left[  t,T\right]  }\leq
R=C_{0}\left(  1+R_{0}\right)
\]
with $C_{0}$ depending only on $M_{0,T}$, $M_{0}$, $L_{0}$ $T$, $\alpha$,
$\beta$ and $\Lambda_{\alpha}\left(  g\right)  .$ Let $\bar{h}=\min\left\{
T-t,1\right\}  .$ Then%
\[
X_{s}^{t,x}=x+\int_{t}^{s}b(r,X_{r}^{t,x})dr+\int_{t}^{s}\sigma(r,X_{r}%
^{t,x})dg\left(  r\right)  \in K\left(  s\right)  ,\quad\forall~s\in\left[
t,t+\bar{h}\right]  .
\]
We clearly have for all $s\in\left[  t,t+\bar{h}\right]  $%
\[
X_{s}^{t,x}=x+\int_{t}^{s}\left[  b(t,x)+U\left(  r\right)  \right]
dr+\int_{t}^{s}\left[  \sigma(t,x)+V\left(  r\right)  \right]  dg\left(
r\right)
\]
where%
\[
U\left(  r\right)  =b(r,X_{r}^{t,x})-b\left(  t,x\right)  \quad\text{and}\quad
V\left(  r\right)  =\sigma(r,X_{r}^{t,x})-\sigma\left(  t,x\right)  .
\]
Clearly $U$ and $V$ satisfy (\ref{UV}).\medskip\newline$\left(  \mathbf{jj}%
\right)  \Rightarrow\left(  \mathbf{jjj}\right)  :$\newline Let $\left\vert
x\right\vert \leq R.$ We verify that
\[
Q\left(  s\right)  =\int_{t}^{s}U\left(  r\right)  dr+\int_{t}^{s}V\left(
r\right)  dg\left(  r\right)
\]
satisfies the Holder conditions from the definition of the contingency
property. Indeed by we have%
\[
\left\vert \int_{t}^{s}U\left(  r\right)  dr\right\vert =\left\vert \int
_{t}^{s}\left[  U\left(  r\right)  -U\left(  t\right)  \right]  dr\right\vert
\leq D_{R}~\left\vert s-t\right\vert ^{2-\alpha}%
\]
and%
\begin{align*}
\left\vert \int_{t}^{s}V\left(  \theta\right)  dg\left(  \theta\right)
\right\vert  & \leq\left\vert \int_{t}^{s}\left[  V\left(  \theta\right)
-V\left(  t\right)  \right]  dg\left(  \theta\right)  \right\vert \\
& \leq\Lambda_{\alpha}(g)%
%TCIMACRO{\dint _{t}^{s}}%
%BeginExpansion
{\displaystyle\int_{t}^{s}}
%EndExpansion
\left(  \dfrac{|V\left(  \theta\right)  -V\left(  t\right)  |}{(\theta
-t)^{\alpha}}+\alpha%
%TCIMACRO{\dint _{t}^{\theta}}%
%BeginExpansion
{\displaystyle\int_{t}^{\theta}}
%EndExpansion
\dfrac{|V(\theta)-V(u)|}{(\theta-u)^{\alpha+1}}du\right)  d\theta\\
& \leq C_{R}^{\left(  1\right)  }\left(  s-t\right)  ^{1+\min\left\{
\beta-\alpha,1-2\alpha\right\}  }%
\end{align*}
Hence%
\[
\left\vert Q\left(  s\right)  \right\vert \leq C_{R}^{\left(  2\right)
}\left(  s-t\right)  ^{1+\min\left\{  \beta-\alpha,1-2\alpha\right\}  }.
\]
Now we prove that $Q$ is $\left(  1-\alpha\right)  -$Holder continuous on
$\left[  t,t+\bar{h}\right]  .$ Let $t\leq\tau\leq s\leq t+\bar{h}.$ We have%
\begin{align*}
\left\vert \int_{\tau}^{s}V\left(  \theta\right)  dg\left(  \theta\right)
\right\vert  & \leq\left\vert \int_{\tau}^{s}\left[  V\left(  \theta\right)
-V\left(  t\right)  \right]  dg\left(  \theta\right)  \right\vert \\
& \leq\Lambda_{\alpha}(g)%
%TCIMACRO{\dint _{\tau}^{s}}%
%BeginExpansion
{\displaystyle\int_{\tau}^{s}}
%EndExpansion
\left(  \dfrac{|V\left(  \theta\right)  -V\left(  t\right)  |}{(\theta
-\tau)^{\alpha}}+\alpha%
%TCIMACRO{\dint _{\tau}^{\theta}}%
%BeginExpansion
{\displaystyle\int_{\tau}^{\theta}}
%EndExpansion
\dfrac{|V(\theta)-V(u)|}{(\theta-u)^{\alpha+1}}du\right)  d\theta\\
& \leq C_{R}^{\left(  3\right)  }\left(  s-\tau\right)  ^{1-\alpha}%
\end{align*}
and therefore
\begin{align*}
\left\vert Q\left(  s\right)  -Q\left(  \tau\right)  \right\vert  &
\leq\left\vert \int_{\tau}^{s}U\left(  r\right)  dr\right\vert +\left\vert
\int_{\tau}^{s}V\left(  \theta\right)  dg\left(  \theta\right)  \right\vert
\\
& \leq C_{R}^{\left(  4\right)  }\left(  s-\tau\right)  ^{1-\alpha}.
\end{align*}
$\left(  \mathbf{jjj}\right)  \Rightarrow\left(  \mathbf{j}\right)  :$\newline

Let us fix $t\in\left[  0,T\right]  $, $x\in K\left(  t\right)  $ and
$0<\varepsilon\leq1.$ Let $R_{0}>0$ be such that $\left\vert x\right\vert \leq
R_{0}.$

We denote by $\mathcal{A}_{\varepsilon}\left(  t,x\right)  $ the set of pairs
$\left(  T_{X},~X\right)  $ where $T_{X}\in\left[  0,T\right]  $ and
$X:\left[  t,T_{X}\right]  \rightarrow\mathbb{R}^{d}$ is a Holder continuous
function satisfying\medskip\newline$\left(  1\right)  \quad X_{t}=x,$
$X_{s}\in K\left(  s\right)  $ for all $s\in\left[  t,T_{X}\right]  ,$ and
there exists a positive constant $B_{0}\geq R_{0}$\ \text{depending only
on\newline }$R_{0}$, $L_{R_{0}}$, $M_{0,T}$, $M_{0}$, $L_{0}$ $T$, $\alpha$,
$\beta$ and $\Lambda_{\alpha}\left(  g\right)  ,$ such that%
\[
\left\Vert X\right\Vert _{1-\alpha;\left[  t,T_{X}\right]  }\leq B_{0}%
\]
$\left(  2\right)  \quad$The error function $\xi:\left[  t,T_{X}\right]
\rightarrow\mathbb{R}^{d}$%
\[
\xi\left(  s\right)  =X_{s}-x-%
%TCIMACRO{\dint _{t}^{s}}%
%BeginExpansion
{\displaystyle\int_{t}^{s}}
%EndExpansion
b\left(  r,X_{r}\right)  dr-%
%TCIMACRO{\dint _{t}^{s}}%
%BeginExpansion
{\displaystyle\int_{t}^{s}}
%EndExpansion
\sigma\left(  r,X_{r}\right)  dg\left(  r\right)  ,\quad s\in\left[
t,T_{X}\right]  ,
\]
satisfies%
\[%
\begin{array}
[c]{ll}%
\left(  a\right)  \quad & \left\vert \xi\left(  s\right)  \right\vert
\leq\varepsilon\left(  s-t\right)  ,\ \text{for\ all}\ s\in\left[
t,T_{X}\right]  ,\medskip\\
\left(  b\right)  \quad & \left\vert \xi\left(  \tau\right)  -\xi\left(
s\right)  \right\vert \leq D_{0}\left\vert \tau-s\right\vert ^{1-\alpha}%
,\quad\text{for all }s,\tau\in\left[  t,T_{X}\right]  ,
\end{array}
\]
where the constant $D_{0}$ depends \text{ only on }$R_{0},L_{R_{0}}$,
$M_{0,T}$, $M_{0}$, $L_{0}$ $T$, $\alpha$, $\beta$ and $\Lambda_{\alpha
}\left(  g\right)  .$

Remark that $B_{0}$ and $D_{0}$ are independent of $\varepsilon$.$\smallskip
$\newline The set $\mathcal{A}_{\varepsilon}\left(  t,x\right)  $ is not empty
because we can find $\left(  t,x\right)  \in\mathcal{A}_{\varepsilon}\left(
t,x\right)  $ .\newline$\mathcal{A}_{\varepsilon}\left(  t,x\right)  $ is an
inductive set for the order relation%
\[
\left(  T_{X_{1}},X_{1}\left(  \cdot\right)  \right)  \preceq\left(  T_{X_{2}%
},X_{2}\left(  \cdot\right)  \right)
\]
defined by%
\[
T_{X_{1}}\leq T_{X_{2}}\text{ and }X_{2}\left(  \cdot\right)  \left\vert
_{\left[  t,T_{X_{1}}\right]  }\right.  =X_{1}\left(  \cdot\right)  .
\]
Zorn's Lemma implies that there exists a maximal element $\left(
T^{\varepsilon},X^{\varepsilon}\right)  $ $\in\mathcal{A}_{\varepsilon}\left(
t,x\right)  $ . We shall prove by \textit{reductio ad absurdum} that
$T^{\varepsilon}=T.$

Assume that $T^{\varepsilon}<T.$ Denote $X_{T^{\varepsilon}}^{\varepsilon
}=x^{\varepsilon}.$ We have $\left\Vert X^{\varepsilon}\right\Vert
_{1-\alpha;\left[  t,T_{X}\right]  }\leq B_{0}$ and in particular
\[
\left\vert x^{\varepsilon}\right\vert \leq B_{0}.
\]
We know from the hypotheses that $\left(  b\left(  T^{\varepsilon
},x^{\varepsilon}\right)  \text{, }\sigma\left(  T^{\varepsilon}%
,x^{\varepsilon}\right)  \right)  $ is $\left(  1-\alpha\right)  -$fractional
$g-$contingent to $K$ in $\left(  T^{\varepsilon},x^{\varepsilon}\right)  $ ,
i.e. there exist $\bar{h}_{\varepsilon}>0$ sufficiently small (for moment
$0<\bar{h}_{\varepsilon}<\left(  T-T^{\varepsilon}\right)  \wedge1$),
$Q^{\varepsilon}:\left[  T^{\varepsilon},T^{\varepsilon}+\bar{h}_{\varepsilon
}\right]  \rightarrow%
%TCIMACRO{\U{211d} }%
%BeginExpansion
\mathbb{R}
%EndExpansion
^{d},$ two constants $G_{0}=G_{B_{0}}$, $\tilde{G}_{0}=\tilde{G}_{B_{0}}>0$
independent of $\left(  T^{\varepsilon},\bar{h}_{\varepsilon}\right)  $ and a
constant $\gamma=\gamma_{B_{0}}\in\left(  0,1\right)  $ also independent of
$\left(  T^{\varepsilon},\bar{h}_{\varepsilon}\right)  $ (the constants
$G_{0}$,$\tilde{G}_{0}$, $\gamma$ depend only on $R_{0},L_{R_{0}},M_{0,T}$,
$M_{0}$, $L_{0}$ $T$, $\alpha$, $\beta$) such that for all $s,\tau\in\left[
T^{\varepsilon},T^{\varepsilon}+\bar{h}_{\varepsilon}\right]  $%
\[
\left\vert Q^{\varepsilon}\left(  \tau\right)  -Q^{\varepsilon}\left(
s\right)  \right\vert \leq G_{0}\left\vert \tau-s\right\vert ^{1-\alpha}%
\quad\text{and}\quad\left\vert Q^{\varepsilon}\left(  s\right)  \right\vert
\leq\tilde{G}_{0}\left\vert s-T^{\varepsilon}\right\vert ^{1+\gamma}%
\]
and satisfying for all $s\in\left[  T^{\varepsilon},T^{\varepsilon}+\bar
{h}_{\varepsilon}\right]  $%
\[
x^{\varepsilon}+\left(  s-T^{\varepsilon}\right)  b\left(  T^{\varepsilon
},x^{\varepsilon}\right)  +\sigma\left(  T^{\varepsilon},x^{\varepsilon
}\right)  \left[  g\left(  s\right)  -g\left(  T^{\varepsilon}\right)
\right]  +Q^{\varepsilon}\left(  s\right)  \in K\left(  s\right)  ,
\]
We set $S^{\varepsilon}=T^{\varepsilon}+\bar{h}_{\varepsilon}$ and we define
$\hat{X}^{\varepsilon}:\left[  t,S^{\varepsilon}\right]  \rightarrow K$ \ as a
extension of $X^{\varepsilon}$ by
\[
\hat{X}^{\varepsilon}\left(  s\right)  =\left\{
\begin{array}
[c]{l}%
X^{\varepsilon}\left(  s\right)  ,\quad\text{if }s\in\left[  t,T^{\varepsilon
}\right]  ,\medskip\\
x^{\varepsilon}+\left(  s-T^{\varepsilon}\right)  b\left(  T^{\varepsilon
},x^{\varepsilon}\right)  +\sigma\left(  T^{\varepsilon},x^{\varepsilon
}\right)  \left(  g\left(  s\right)  -g\left(  T^{\varepsilon}\right)
\right)  +Q^{\varepsilon}\left(  s\right)  ,\medskip\\
\qquad\qquad\qquad\qquad\qquad\qquad\qquad\qquad\qquad\quad\text{if }%
s\in\left[  T^{\varepsilon},S^{\varepsilon}\right]  .
\end{array}
\right.
\]
We will prove that the extension $\left(  S^{\varepsilon},\hat{X}%
^{\varepsilon}\right)  \in\mathcal{A}_{\varepsilon}\left(  t,x\right)
$.\newline\textit{Step 1: }Clearly $\hat{X}_{t}^{\varepsilon}=x$ and $\hat
{X}_{s}^{\varepsilon}\in K\left(  s\right)  $ for all $s\in\left[  t,T\right]
.$

Let us show that $\left\Vert \hat{X}^{\varepsilon}\right\Vert _{1-\alpha
;\left[  t,S^{\varepsilon}\right]  }\leq B_{0}^{\left(  1\right)  }$ where
$B_{0}^{\left(  1\right)  }\geq R_{0}$ and \ $B_{0}^{\left(  1\right)  }$
\text{depends only on }$R_{0},L_{R_{0}}$, $M_{0,T}$, $M_{0}$, $L_{0}$ $T$,
$\alpha$, $\beta$ and $\Lambda_{\alpha}\left(  g\right)  .$

Let $T^{\varepsilon}\leq s\leq\tau\leq S^{\varepsilon}.$ Then
\begin{equation}%
\begin{array}
[c]{r}%
\left\vert \hat{X}_{\tau}^{\varepsilon}-\hat{X}_{s}^{\varepsilon}\right\vert
\leq\left\vert \tau-s\right\vert \left\vert b\left(  T^{\varepsilon
},x^{\varepsilon}\right)  \right\vert +\Lambda_{\alpha}\left(  g\right)
\Gamma\left(  \alpha\right)  \left\vert \sigma\left(  T^{\varepsilon
},x^{\varepsilon}\right)  \right\vert \left\vert \tau-s\right\vert ^{1-\alpha
}+\left\vert Q^{\varepsilon}\left(  \tau\right)  -Q^{\varepsilon}\left(
s\right)  \right\vert \medskip\\
\multicolumn{1}{l}{\quad\leq\left\vert \tau-s\right\vert L_{0}\left(
1+\left\vert x^{\varepsilon}\right\vert \right)  +\dfrac{\Lambda_{\alpha
}\left(  g\right)  }{\alpha}M_{0,T}\left(  1+\left\vert x^{\varepsilon
}\right\vert \right)  \left\vert \tau-s\right\vert ^{1-\alpha}+G_{R}\left\vert
\tau-s\right\vert ^{1-\alpha}\medskip}\\
\multicolumn{1}{l}{\quad\leq C_{R_{0}}\left\vert \tau-s\right\vert ^{1-\alpha
},}%
\end{array}
\label{holder}%
\end{equation}
and $\left\vert \hat{X}_{\tau}^{\varepsilon}\right\vert \leq\left\vert \hat
{X}_{\tau}^{\varepsilon}-\hat{X}_{T^{\varepsilon}}^{\varepsilon}\right\vert
+\left\vert x^{\varepsilon}\right\vert \leq C_{R_{0}}T^{1-\alpha}+R_{0}.$
Hence%
\begin{align*}
\left\Vert \hat{X}^{\varepsilon}\right\Vert _{1-\alpha;\left[
t,S^{\varepsilon}\right]  }  &  \leq\left\Vert \hat{X}^{\varepsilon
}\right\Vert _{1-\alpha;\left[  t,T^{\varepsilon}\right]  }+\left\Vert \hat
{X}^{\varepsilon}\right\Vert _{1-\alpha;\left[  T^{\varepsilon},S^{\varepsilon
}\right]  }\\
&  \leq B_{0}+C_{R_{0}}T^{1-\alpha}+R_{0}+C_{R_{0}}\overset{def}{=}%
B_{0}^{\left(  1\right)  }~.
\end{align*}
\textit{Step 2: the error function.}

Let the error functions $\xi^{\varepsilon}:\left[  t,T^{\varepsilon}\right]
\rightarrow\mathbb{R}^{d}$ and $\hat{\xi}^{\varepsilon}:\left[
t,S^{\varepsilon}\right]  \rightarrow\mathbb{R}^{d},$%
\begin{align*}
\xi^{\varepsilon}\left(  s\right)   &  =X_{s}^{\varepsilon}-x-%
%TCIMACRO{\dint _{t}^{s}}%
%BeginExpansion
{\displaystyle\int_{t}^{s}}
%EndExpansion
b\left(  r,X_{r}^{\varepsilon}\right)  dr-%
%TCIMACRO{\dint _{t}^{s}}%
%BeginExpansion
{\displaystyle\int_{t}^{s}}
%EndExpansion
\sigma\left(  r,X_{r}^{\varepsilon}\right)  dg\left(  r\right)  .\\
\hat{\xi}^{\varepsilon}\left(  s\right)   &  =\hat{X}_{s}^{\varepsilon}-x-%
%TCIMACRO{\dint _{t}^{s}}%
%BeginExpansion
{\displaystyle\int_{t}^{s}}
%EndExpansion
b\left(  r,\hat{X}_{r}^{\varepsilon}\right)  dr-%
%TCIMACRO{\dint _{t}^{s}}%
%BeginExpansion
{\displaystyle\int_{t}^{s}}
%EndExpansion
\sigma\left(  r,\hat{X}_{r}^{\varepsilon}\right)  dg\left(  r\right)  .
\end{align*}
Clearly $\left\vert \hat{\xi}^{\varepsilon}\left(  s\right)  \right\vert
=\left\vert \xi^{\varepsilon}\left(  s\right)  \right\vert \leq\varepsilon
\left(  s-t\right)  $ for all $s\in\left[  t,T^{\varepsilon}\right]  .$

Let $s\in\left[  T^{\varepsilon},S^{\varepsilon}\right]  .$ Using Lemma
\ref{aux} (the inequalities (\ref{aux1}) with $Y_{r}=\hat{X}_{r}^{\varepsilon
}$ and $t=T^{\varepsilon},$ we have\medskip\newline$\left\vert \hat{\xi
}^{\varepsilon}\left(  s\right)  \right\vert \leq\left\vert x^{\varepsilon}-x-%
%TCIMACRO{\dint _{t}^{T^{\varepsilon}}}%
%BeginExpansion
{\displaystyle\int_{t}^{T^{\varepsilon}}}
%EndExpansion
b\left(  r,X_{r}^{\varepsilon}\right)  dr-%
%TCIMACRO{\dint _{t}^{T^{\varepsilon}}}%
%BeginExpansion
{\displaystyle\int_{t}^{T^{\varepsilon}}}
%EndExpansion
\sigma\left(  r,X_{r}^{\varepsilon}\right)  dg\left(  r\right)  \right\vert
\medskip\newline\ +\left\vert \left(  s-T^{\varepsilon}\right)  b\left(
T^{\varepsilon},x^{\varepsilon}\right)  +\sigma\left(  T^{\varepsilon
},x^{\varepsilon}\right)  \left(  g\left(  s\right)  -g\left(  T^{\varepsilon
}\right)  \right)  -%
%TCIMACRO{\dint _{T^{\varepsilon}}^{s}}%
%BeginExpansion
{\displaystyle\int_{T^{\varepsilon}}^{s}}
%EndExpansion
b\left(  r,\hat{X}_{r}^{\varepsilon}\right)  dr-%
%TCIMACRO{\dint _{T^{\varepsilon}}^{s}}%
%BeginExpansion
{\displaystyle\int_{T^{\varepsilon}}^{s}}
%EndExpansion
\sigma\left(  r,\hat{X}_{r}^{\varepsilon}\right)  dg\left(  r\right)
\right\vert \medskip\newline\ +\left\vert Q^{\varepsilon}\left(  s\right)
\right\vert \medskip\newline\leq\varepsilon\left(  T^{\varepsilon}-t\right)
+\left\vert
%TCIMACRO{\dint _{T^{\varepsilon}}^{s}}%
%BeginExpansion
{\displaystyle\int_{T^{\varepsilon}}^{s}}
%EndExpansion
\left[  b\left(  T^{\varepsilon},x^{\varepsilon}\right)  -b\left(  r,\hat
{X}_{r}^{\varepsilon}\right)  \right]  dr+%
%TCIMACRO{\dint _{T^{\varepsilon}}^{s}}%
%BeginExpansion
{\displaystyle\int_{T^{\varepsilon}}^{s}}
%EndExpansion
\left[  \sigma\left(  T^{\varepsilon},x^{\varepsilon}\right)  -\sigma\left(
r,\hat{X}_{r}^{\varepsilon}\right)  \right]  dg\left(  r\right)  \right\vert
\medskip\newline%
\begin{array}
[c]{c}%
\quad
\end{array}
+\tilde{G}_{0}\left(  s-T^{\varepsilon}\right)  ^{1+\gamma}\medskip
\newline\leq\varepsilon\left(  T^{\varepsilon}-t\right)  +C_{R_{0}}^{\left(
1\right)  }\left(  s-T^{\varepsilon}\right)  ^{2-\alpha}+C_{R_{0}}^{\left(
2\right)  }\left(  s-T^{\varepsilon}\right)  ^{1+\min\left\{  \beta
-\alpha,1-2\alpha\right\}  }+\tilde{G}_{0}\left(  s-T^{\varepsilon}\right)
^{1+\gamma}\medskip\newline\leq\varepsilon\left(  T^{\varepsilon}-t\right)
+\varepsilon\left(  s-T^{\varepsilon}\right)  =\varepsilon\left(  s-t\right)
\medskip\newline$for $\bar{h}_{\varepsilon}$ sufficiently small such that%
\[
C_{R_{0}}^{\left(  1\right)  }\bar{h}_{\varepsilon}^{1-\alpha}+C_{R_{0}%
}^{\left(  2\right)  }\bar{h}_{\varepsilon}^{\min\left\{  \beta-\alpha
,1-2\alpha\right\}  }+\tilde{G}_{0}\bar{h}_{\varepsilon}^{\gamma}%
\leq\varepsilon.
\]
Hence%
\[
\left\vert \hat{\xi}^{\varepsilon}\left(  s\right)  \right\vert \leq
\varepsilon\left(  s-t\right)  \text{ for all }s\in\left[  t,S^{\varepsilon
}\right]  .
\]
Let now $T^{\varepsilon}\leq\tau<s\leq S^{\varepsilon}.$ Then by Lemma
\ref{aux} (the inequalities (\ref{aux2}) with $Y_{r}=\hat{X}_{r}^{\varepsilon
}~,$ $t=T^{\varepsilon}$ and $x^{\varepsilon}=\hat{X}_{T^{\varepsilon}%
}^{\varepsilon}$ we have\medskip\newline$\left\vert \hat{\xi}^{\varepsilon
}\left(  s\right)  -\hat{\xi}^{\varepsilon}\left(  \tau\right)  \right\vert
$\medskip\newline$%
\begin{array}
[c]{c}%
\quad
\end{array}
\leq\left\vert
%TCIMACRO{\dint _{\tau}^{s}}%
%BeginExpansion
{\displaystyle\int_{\tau}^{s}}
%EndExpansion
\left[  b\left(  T^{\varepsilon},x^{\varepsilon}\right)  -b\left(  r,\hat
{X}_{r}^{\varepsilon}\right)  \right]  dr+%
%TCIMACRO{\dint _{\tau}^{s}}%
%BeginExpansion
{\displaystyle\int_{\tau}^{s}}
%EndExpansion
\left[  \sigma\left(  T^{\varepsilon},x^{\varepsilon}\right)  -\sigma\left(
r,\hat{X}_{r}^{\varepsilon}\right)  \right]  dg_{r}\right\vert +\left\vert
Q^{\varepsilon}\left(  s\right)  -Q^{\varepsilon}\left(  \tau\right)
\right\vert $\medskip\newline$%
\begin{array}
[c]{c}%
\quad
\end{array}
\leq C_{R_{0}}^{\left(  3\right)  }\left(  s-\tau\right)  +C_{R_{0}}^{\left(
4\right)  }\left(  s-\tau\right)  ^{1-\alpha}+G_{0}\left\vert \tau
-s\right\vert ^{1-\alpha}$\medskip\newline$%
\begin{array}
[c]{c}%
\quad
\end{array}
\leq C_{R_{0}}\left(  s-\tau\right)  ^{1-\alpha}.$\medskip\newline From the
definition of $\mathcal{A}_{\varepsilon}\left(  t,x\right)  ,$ for $\tau
,s\in\left[  t,T^{\varepsilon}\right]  $
\[
\left\vert \hat{\xi}^{\varepsilon}\left(  s\right)  -\hat{\xi}^{\varepsilon
}\left(  \tau\right)  \right\vert =\left\vert \xi^{\varepsilon}\left(
s\right)  -\xi^{\varepsilon}\left(  \tau\right)  \right\vert \leq
D_{0}~\left(  s-\tau\right)  ^{1-\alpha},\;\;\forall~\tau,s\in\left[
t,T^{\varepsilon}\right]  .
\]
We conclude%
\[
\left\vert \hat{\xi}^{\varepsilon}\left(  s\right)  -\hat{\xi}^{\varepsilon
}\left(  \tau\right)  \right\vert \leq\left(  C_{R_{0}}\vee D_{0}\right)
\left\vert \tau-s\right\vert ^{1-\alpha},\;\;\forall~\tau,s\in\left[
t,S^{\varepsilon}\right]  .
\]
We arrived to prove that $\left(  S^{\varepsilon},\hat{X}^{\varepsilon
}\right)  $ is proper extension of $\left(  T^{\varepsilon},X^{\varepsilon
}\right)  ,$ that contradicts the maximallity of $\left(  T^{\varepsilon
},X^{\varepsilon}\right)  $ in $\mathcal{A}_{\varepsilon}\left(  t,x\right)
.$ Therefore $T^{\varepsilon}=T.$

Let $\left(  T,~X^{\varepsilon}\right)  $ be a maximal element of
$\mathcal{A}_{\varepsilon}\left(  t,x\right)  .$ Then from the definition of
$\mathcal{A}_{\varepsilon}\left(  t,x\right)  $ we have $X_{t}^{\varepsilon
}=x,$ $X_{s}^{\varepsilon}\in K\left(  s\right)  $ for all $s\in\left[
t,T\right]  ,$ and there exists a positive constant $B_{0}\geq R_{0}$
\ \text{depending only on }$R_{0},L_{R_{0}}$, $M_{0,T}$, $M_{0}$, $L_{0}$ $T$,
$\alpha$, $\beta$ and $\Lambda_{\alpha}\left(  g\right)  ,$ such that%
\[
\left\Vert X^{\varepsilon}\right\Vert _{1-\alpha;\left[  t,T\right]  }\leq
B_{0}\,.
\]
The error function $\xi^{\varepsilon}:\left[  t,T\right]  \rightarrow
\mathbb{R}^{d},\mathbb{\ }$%
\begin{equation}
\xi^{\varepsilon}\left(  s\right)  =X_{s}^{\varepsilon}-x-%
%TCIMACRO{\dint _{t}^{s}}%
%BeginExpansion
{\displaystyle\int_{t}^{s}}
%EndExpansion
b\left(  r,X_{r}^{\varepsilon}\right)  dr-%
%TCIMACRO{\dint _{t}^{s}}%
%BeginExpansion
{\displaystyle\int_{t}^{s}}
%EndExpansion
\sigma\left(  r,X_{r}^{\varepsilon}\right)  dg\left(  r\right)
\label{ec-aprox}%
\end{equation}
satisfies
\[%
\begin{array}
[c]{ll}%
\left(  a\right)  \quad & \left\vert \xi^{\varepsilon}\left(  s\right)
\right\vert \leq\varepsilon\left(  s-t\right)  ,\ \text{for\ all}\ s\in\left[
t,T\right]  ,\medskip\\
\left(  b\right)  \quad & \left\vert \xi^{\varepsilon}\left(  \tau\right)
-\xi^{\varepsilon}\left(  s\right)  \right\vert \leq D_{0}\left\vert
\tau-s\right\vert ^{1-\alpha},\quad\text{for all }s,\tau\in\left[  t,T\right]
,
\end{array}
\]
where the constant $D_{0}$ depends \text{ only on }$R_{0},L_{R_{0}}$,
$M_{0,T}$, $M_{0}$, $L_{0}$ $T$, $\alpha$, $\beta$ and $\Lambda_{\alpha
}\left(  g\right)  .$

We now estimate $\left\Vert \xi^{\varepsilon}\right\Vert _{\alpha
,\infty;\left[  t,T\right]  }~.$ Let $\lambda\geq0.$ Then\medskip
\newline$\left\vert \!\left\vert \!\left\vert \xi^{\varepsilon}\right\vert
\!\right\vert \!\right\vert _{\alpha,\lambda;\left[  t,T\right]  }%
\leq\left\Vert \xi^{\varepsilon}\right\Vert _{\alpha,\infty;\left[
t,T\right]  }$\medskip\newline$%
\begin{array}
[c]{c}%
\quad\quad\quad\quad
\end{array}
=\underset{s\in\left[  t,T\right]  }{\sup}\left\{  \left\vert \xi
^{\varepsilon}\left(  s\right)  \right\vert +%
%TCIMACRO{\dint _{t}^{s}}%
%BeginExpansion
{\displaystyle\int_{t}^{s}}
%EndExpansion
\dfrac{\left\vert \xi^{\varepsilon}\left(  s\right)  -\xi^{\varepsilon}\left(
r\right)  \right\vert }{\left(  s-r\right)  ^{\alpha+1}}dr\right\}  $%
\medskip\newline$%
\begin{array}
[c]{c}%
\quad\quad\quad\quad
\end{array}
\leq\varepsilon\left(  T-t\right)  +\underset{s\in\left[  t,T\right]  }{\sup}%
%TCIMACRO{\dint _{t}^{s}}%
%BeginExpansion
{\displaystyle\int_{t}^{s}}
%EndExpansion
\left\vert \xi_{s}^{\varepsilon}-\xi_{r}^{\varepsilon}\right\vert ^{\frac
{1}{2}-\alpha}\dfrac{\left\vert \xi_{s}^{\varepsilon}-\xi_{r}^{\varepsilon
}\right\vert ^{\frac{1}{2}+\alpha}}{\left(  s-r\right)  ^{\alpha+1}}%
dr$\medskip\newline$%
\begin{array}
[c]{c}%
\quad\quad\quad\quad
\end{array}
\leq\varepsilon\left(  T-t\right)  +\left[  2\varepsilon\left(  T-t\right)
\right]  ^{\frac{1}{2}-\alpha}\underset{s\in\left[  t,T\right]  }{\sup}%
%TCIMACRO{\dint _{t}^{s}}%
%BeginExpansion
{\displaystyle\int_{t}^{s}}
%EndExpansion
\dfrac{D_{0}^{\frac{1}{2}+\alpha}\left(  s-r\right)  ^{\left(  1-\alpha
\right)  \left(  \frac{1}{2}+\alpha\right)  }}{\left(  s-r\right)  ^{1+\alpha
}}dr$\medskip\newline$%
\begin{array}
[c]{c}%
\quad\quad\quad\quad
\end{array}
=\varepsilon\left(  T-t\right)  +\left[  2\varepsilon\left(  T-t\right)
\right]  ^{\frac{1}{2}-\alpha}\frac{D_{0}^{\frac{1}{2}+\alpha}}{\left(
\frac{1}{2}-\alpha\right)  \left(  1+\alpha\right)  }\left(  T-t\right)
^{\left(  \frac{1}{2}-\alpha\right)  \left(  1+\alpha\right)  }$%
\medskip\newline$%
\begin{array}
[c]{c}%
\quad\quad\quad\quad
\end{array}
\leq C_{R_{0}}~\varepsilon^{\frac{1}{2}-\alpha}~.$\medskip\newline It remains
now to prove that the limit of the sequence $X^{\varepsilon}$ exists as
$\varepsilon\rightarrow0$ and this limit is a solution to the differential
equation (\ref{det}) \newline Let $0<\varepsilon,\eta\leq1$. Using the
estimates (\ref{est-l1-2}) and (\ref{est-b2}), we get\medskip\newline%
$\left\vert \!\left\vert \!\left\vert X^{\varepsilon}-X^{\eta}\right\vert
\!\right\vert \!\right\vert _{\alpha,\lambda;\left[  t,T\right]  }$%
\medskip\newline$\leq\left\vert \!\left\vert \!\left\vert F_{t,\cdot}^{\left(
b\right)  }\left(  X^{\varepsilon}\right)  -F_{t,\cdot}^{\left(  b\right)
}\left(  X^{\eta}\right)  \right\vert \!\right\vert \!\right\vert
_{\alpha,\lambda;\left[  t,T\right]  }+\left\vert \!\left\vert \!\left\vert
G_{t,\cdot}^{\left(  \sigma\right)  }\left(  X^{\varepsilon}\right)
-G_{t,\cdot}^{\left(  \sigma\right)  }\left(  X^{\eta}\right)  \right\vert
\!\right\vert \!\right\vert _{\alpha,\lambda;\left[  t,T\right]  }+\left\vert
\!\left\vert \!\left\vert \xi^{\varepsilon}-\xi^{\eta}\right\vert
\!\right\vert \!\right\vert _{\alpha,\lambda;\left[  t,T\right]  }$%
\medskip\newline$\leq\dfrac{C_{R}^{(b3)}}{\lambda^{\alpha}}\left\vert
\!\left\vert \!\left\vert X^{\varepsilon}-X^{\eta}\right\vert \!\right\vert
\!\right\vert _{\alpha,\lambda;\left[  t,T\right]  }+\dfrac{C_{R}^{(\sigma
3)}\Lambda_{\alpha}(g)}{\lambda^{1-2\alpha}}\left(  1+\Delta_{\left[
t,T\right]  }\left(  X^{\varepsilon}\right)  +\Delta_{\left[  t,T\right]
}\left(  X^{\eta}\right)  \right)  \left\vert \!\left\vert \!\left\vert
X^{\varepsilon}-X^{\eta}\right\vert \!\right\vert \!\right\vert _{\alpha
,\lambda;\left[  t,T\right]  }\medskip\newline%
\begin{array}
[c]{c}%
\quad\quad\quad
\end{array}
+C_{R_{0}}~\varepsilon^{\frac{1}{2}-\alpha}+C_{R_{0}}~\eta^{\frac{1}{2}%
-\alpha}\medskip$\newline

Since $\left\Vert X^{\varepsilon}\right\Vert _{1-\alpha;\left[  t,T\right]
}\leq B_{0}~,\,0<\alpha<\frac{\delta}{1+\delta},$ then%
\[%
\begin{array}
[c]{ll}%
\Delta_{\left[  t,T\right]  }\left(  X^{\varepsilon}\right)  =\sup
_{r\in\left[  t,T\right]  }%
%TCIMACRO{\dint _{t}^{r}}%
%BeginExpansion
{\displaystyle\int_{t}^{r}}
%EndExpansion
\dfrac{\left\vert X_{r}^{\varepsilon}-X_{s}^{\varepsilon}\right\vert ^{\delta
}}{\left(  r-s\right)  ^{1+\alpha}} & \leq B_{0}\sup_{r\in\left[  t,T\right]
}%
%TCIMACRO{\dint _{t}^{r}}%
%BeginExpansion
{\displaystyle\int_{t}^{r}}
%EndExpansion
\dfrac{(r-s)^{\left(  1-\alpha\right)  \delta}}{(r-s)^{\alpha+1}}ds\medskip\\
& \leq B_{0}\dfrac{T^{\delta-\alpha\left(  1+\delta\right)  }}{\delta
-\alpha\left(  1+\delta\right)  }\medskip\\
& \leq\dfrac{B_{0}\left(  1+T\right)  }{\delta-\alpha\left(  1+\delta\right)
}%
\end{array}
\]
and therefore%
\begin{align*}
\left\vert \!\left\vert \!\left\vert X^{\varepsilon}-X^{\eta}\right\vert
\!\right\vert \!\right\vert _{\alpha,\lambda;\left[  t,T\right]  }  &
\leq\frac{C_{R_{0}}^{\left(  3\right)  }}{\delta-\alpha\left(  1+\delta
\right)  }\left(  \dfrac{1}{\lambda^{\alpha}}+\dfrac{1}{\lambda^{1-2\alpha}%
}\right)  \left\vert \!\left\vert \!\left\vert X^{\varepsilon}-X^{\eta
}\right\vert \!\right\vert \!\right\vert _{\alpha,\lambda;\left[  t,T\right]
}\medskip\\
&  +C_{R_{0}}~\varepsilon^{\frac{1}{2}-\alpha}+C_{R_{0}}~\eta^{\frac{1}%
{2}-\alpha}%
\end{align*}

Let $\lambda=\bar{\lambda}\geq1$ such that%
\[
\frac{C_{R_{0}}^{\left(  3\right)  }}{\delta-\alpha\left(  1+\delta\right)
}\left(  \dfrac{1}{\bar{\lambda}^{\alpha}}+\dfrac{1}{\bar{\lambda}^{1-2\alpha
}}\right)  \leq\frac{1}{2}.
\]
We deduce%
\[
\left\Vert X^{\varepsilon}-X^{\eta}\right\Vert _{\infty;\left[  t,T\right]
}\leq e^{\bar{\lambda}T}\left\vert \!\left\vert \!\left\vert X^{\varepsilon
}-X^{\eta}\right\vert \!\right\vert \!\right\vert _{\alpha,\bar{\lambda
};\left[  t,T\right]  }\leq2C_{R_{0}}e^{\bar{\lambda}T}~\left(  \varepsilon
^{\frac{1}{2}-\alpha}+~\eta^{\frac{1}{2}-\alpha}\right)  .
\]
Hence there exists $X^{t,x}$ such that $X^{\varepsilon}\rightarrow X^{t,x}$ in
$C\left(  \left[  t,T\right]  ;K\right)  $ and $X^{\varepsilon}\rightarrow
X^{t,x}$ in $W^{\alpha,\infty}(t,T;\mathbb{R}^{d})$ as $\varepsilon
\rightarrow0.$ Since for all $s,\tau\in\left[  t,T\right]  :$%
\[
\left\Vert X^{\varepsilon}\right\Vert _{\infty;\left[  t,T\right]  }%
+\frac{\left\vert X_{s}^{\varepsilon}-X_{\tau}^{\varepsilon}\right\vert
}{\left\vert s-\tau\right\vert ^{1-\alpha}}\leq B_{0}%
\]
then passing to limit as $\varepsilon\rightarrow0$ we obtain
\[
\left\Vert X^{t,x}\right\Vert _{1-\alpha;\left[  t,T\right]  }\leq B_{0}.
\]
Since $X_{s}^{\varepsilon}\in K\left(  s\right)  $ for all $s\in\left[
t,T\right]  ,$ clearly follows, as $\varepsilon\rightarrow0,$ that
\[
X_{s}\in K\left(  s\right)  \quad\text{for all}\quad s\in\left[  t,T\right]
.
\]
Passing to limit in (\ref{ec-aprox}) we infer that $X^{t,x}$ is a solution of
the differential equation (\ref{det}) starting at $t$ from $x$ and evolving in
the tube $\left\{  \left(  s,y\right)  :s\in\left[  t,T\right]  ,\ y\in
K\left(  s\right)  \right\}  .$\newline The proof is complete.

\hfill
\end{proof}

\section{Proof of the main result}

Let be fixed a parameter $1/2<H<1$. Consider $B=\left\{  B_{t},t\in\left[
0,T\right]  \right\}  $ be a $\mathbb{R}^{k}$ valued fractional Brownian
motion with parameter $H$ defined in a complete probability space $\left(
\Omega,\mathcal{F},\mathbb{P}\right)  .$ From (\ref{covari})it follows that
\[
\mathbb{E}\left(  \left\vert B_{t}-B_{s}\right\vert ^{2}\right)  =\left\vert
t-s\right\vert ^{2H}%
\]
and as a consequence, for any $p\geq1,$%
\[
\left\Vert B_{t}-B_{s}\right\Vert _{p}=\left(  \mathbb{E}\left(  \left\vert
B_{t}-B_{s}\right\vert ^{p}\right)  \right)  ^{1/p}=c_{p}\left\vert
t-s\right\vert ^{H}.
\]
It is known that the random variable%
\[
G=\frac{1}{\Gamma\left(  1-\alpha\right)  }\underset{t<s<r<T}{\sup}\left\vert
\left(  D_{r-}^{1-\alpha}B_{r-}\right)  \left(  s\right)  \right\vert
\]
has moments of all order. As a consequence, As a consequence, if
$u=\{u_{t},\;t\in\lbrack0,T]\}$ is a stochastic process whose trajectories
belong to the space $W^{\alpha,1}(0,T;\mathbb{R}^{d\times k})$, with
$1-H<\alpha<\frac{1}{2}$, the pathwise integral $\int_{0}^{T}u_{s}dB_{s}$
exists in the sense of Definition \ref{def-int} and we have the estimate
\[
\left\vert \int_{0}^{T}u_{s}dB_{s}\right\vert \leq G~\left\Vert u\right\Vert
_{\alpha,1}\text{. }%
\]
Moreover, if the trajectories of the process $u$ belong to the space
$W^{\alpha,\infty}(0,T;\mathbb{R}^{d\times k})$, then the indefinite integral
$U_{t}=\int_{0}^{t}u_{s}dB_{s}$ $\ $\ is Holder continuous of order $1-\alpha
$, and the estimates from Proposition \ref{p-est} hold.

\begin{proof}
[Proof of the Theorem \ref{main}]Considering the previous observations, the
solution follows directly from the deterministic Theorem \ref{detmain}.

\hfill
\end{proof}

\end{document}